\documentclass[a4paper,fleqn]{cas-sc}
\usepackage{amsthm} 
\usepackage[numbers]{natbib}

\usepackage{soul}
\usepackage{xcolor}
\usepackage{algorithm} 
\usepackage{algpseudocode}
\usepackage{amsfonts,amssymb, bm}
\usepackage{graphicx,color,epstopdf,subcaption}
\usepackage{float}
\floatstyle{plain}
\restylefloat{figure}%
\usepackage{booktabs}
\usepackage{bm}
\usepackage{mathtools}

\newtheorem{myremark}{Remark}[section]

\def\XXint#1#2#3{{\setbox0=\hbox{$#1{#2#3}{\int}$}
    \vcenter{\hbox{$#2#3$}}\kern-.5\wd0}}

\def\tx {\tilde{x}}

\def\bi{{\bf i}}

\def\XXint#1#2#3{{\setbox0=\hbox{$#1{#2#3}{\int}$}	\vcenter{\hbox{$#2#3$}}\kern-.5\wd0}}

\begin{document}
\let\WriteBookmarks\relax
\def\floatpagepagefraction{1}
\def\textpagefraction{.001}
\shorttitle{Floquet-Bloch transform based BIE method}
\shortauthors{W. Lu, K. Shen and R. Zhang}

\title [mode = title]{A boundary integral equation method for wave scattering in periodic structures via the Floquet-Bloch transform}          

\tnotemark[1]


\author[1]{Wangtao Lu}[style=chinese]
\ead{wangtaolu@zju.edu.cn}
\cormark[1]

\cortext[cor1]{Corresponding author}

\affiliation[1]{organization={School of Mathematical Sciences, Zhejiang University},
                city={Hangzhou},
                state={Zhejiang},
                country={China}}

\author[1]{Kuanrong Shen}[style=chinese]
\ead{kuanrongshen@zju.edu.cn}

\author[2]{Ruming Zhang}[style=chinese]
\ead{ruming.zhang@tu-berlin.de}

\affiliation[2]{organization={Institute of Mathematics, TU Berlin},
                city={10623 Berlin},
                state={Berlin},
                country={Germany}}


\begin{abstract}
This paper is concerned with the problem of an acoustic wave scattering in a
locally perturbed periodic structure. As the total wavefield is
non-quasi-periodic, effective truncation techniques are pursued for
high-accuracy numerical solvers. We adopt the Green's function for the
background periodic structure to construct a boundary integral equation (BIE) on
an artificial curve enclosing the perturbation. It serves as a transparent
boundary condition (TBC) to truncate the unbounded domain. We develop efficient
algorithms to compute such background Green's functions based on the
Floquet-Bloch transform and its inverse. Spectrally accurate quadrature rules
are developed to discretize the BIE-based TBC. Effective algorithms based on
leap and pullback procedures are further developed to compute the total
wavefield everywhere in the structure. A number of numerical experiments are
carried out to illustrate the efficiency and accuracy of the new solver. They
exhibit that our method for the non-quasi-periodic problem has a time complexity
that is even comparable to that of a single quasi-periodic problem. 
\end{abstract}




\maketitle

\section{Introduction}
Due to their powerful and flexible abilities in manipulating waves,  periodic
structures have seen wide-ranging applications across acoustics,
electromagnetics (optics), elastodynamics, and so on. These practical needs
highly purse fast and accurate numerical simulations for the underlying
mathematical models. Nevertheless, challenges from geometries of the structures
and radiation behaviors of the scattered waves make rigorous mathematical
theories and effective numerical solvers equally difficult to establish in
practice \cite{baolibook22}. Scattering problems in periodic structures have
thus been a hot research topic in communities such as engineering, scientific
computing and applied mathematics, in the past decades. This paper studies an
acoustic wave scattering in a two-dimensional periodic structure with or
without local perturbations. Throughout the rest of this paper, we always assume
that the periodic structures are sandwiched by two homogeneous media in the
transverse direction that is perpendicular to the periodic directions. 

The fundamental problem of a plane wave scattered by a periodic structure has
been studied intensively among the existing literature. The total wavefield is
quasi-periodic so that the problem can be reduced onto a single unit cell.
Moreover, the quasi-periodicity leads to Rayleigh expansions of the total
wavefield, yielding transparent boundary conditions (TBCs) terminating the
transverse direction. Such a boundary value problem has been proved to be
well-posed at all wave frequencies except a countable number of frequencies
\cite{bao1995mathematical}. Readers are referred to \cite{bonnet94guided} for a
rigorous proof of the existence of ill-posed frequencies, the so-called embedded
eigenmodes or bound states in the continuum (BIC), and also to
\cite{hsu2016bound} and the references therein for the extensive applications of
BICs.

Suppose the periodic structure is piecewise homogeneous and the quasi-periodic
problem is well-posed for a given frequency. Among the existing solvers, the
boundary integral equation (BIE) methods are attractive and competitive for
such a scattering problem, mainly because they formulate the problem on the
boundaries of homogeneous domains, thereby reducing the problem dimensionality
by one. Existing BIE methods can generally be classified into two approaches.
The first approach represents the wavefield using single or double-layer potentials
that rely on quasi-periodic Green's functions \cite{barnett2011new}. 
Radiation conditions at
infinity are automatically satisfied and BIEs are established on bounded
surfaces only. The main difficulty is the evaluation of the complicated
quasi-periodic Green's functions, especially when the frequency lies in the
vicinity of Rayleigh anomalies in the sense that plane waves propagating
horizontally appear in the Rayleigh expansions, and the quasi-periodic Green's
functions are not well defined! Readers are referred to \cite{bruno20,zhang2018fft} for effective
algorithms for evaluating the quasi-periodic Green's functions and for more
details of the difficulty. 
The second approach uses free-space Green's
functions, simpler and ready to use, to establish BIEs. As these Green's
functions are not quasi-periodic, BIEs are established on the boundaries of
homogeneities in a unit cell, which always contain irregular corners/edges
\cite{lu2012high,wu2011boundary}. Therefore, the main difficulty comes from the discretizations of
the BIEs. Fortunately, high-accuracy techniques for discretizing BIEs on
piecewise smooth curves have been well-established \cite{alp99, colkre13} so
that the second approach becomes arguably more promising at least for 2D
problems. Readers are referred to \cite{baolibook22} for various numerical
solvers for such  quasi-periodic problems. 

When the periodic structure contains local perturbations or the incident wave is
not quasi-periodic, the total wavefield is not quasi-periodic anymore. It
becomes a fundamental barrier against an effective numerical solver to
accurately truncate the unbounded domain. With perfectly matched layers (PMLs)
\cite{ber94} truncating the transverse direction \cite{chamon09} in advance, existing
truncation techniques for periodic structures in closed waveguides including
TBCs established by Ricatti equations \cite{joly2006exact, yu2022pml},
recursively doubling procedures
\cite{yuan2007recursive,ehrhardt2008numerical,sun2008numerical}, 
or Bloch-mode expansions \cite{helfert1998efficient,hu2009efficient},
become applicable. We point out, in \cite{yu2022pml}, we firstly developed a PML-based
BIE solver for wave scattered by a locally perturbed periodic curve, which can
attain a high-order accuracy; meanwhile, we have found in \cite{zhang2022exponential,zhang2022fast} that
the PML convergence rate may degrade from exponential to algebraic at or close
to half-integer wavenumbers (for $2\pi$-period structures), thereby reducing its overall efficiency. Probably, the
most effective approach should be attributed to the method of Floquet-Bloch (FB)
transform. The wavefield is transformed as an integral of quasi-periodic waves
\cite{lechleiter2017floquet}, each of which can be computed effectively via the
aforementioned solvers. Nevertheless, if the periodic structure contains a local
perturbation, a nonlinear map making the structure periodic again has to be sought
first \cite{lechleiterzhang2017floquet}
to make the FB transform method
applicable. The map complicates the governing equations and significantly
undermines the efficiency of the FB transform, not to mention that such a map may
not even exist if the perturbation destroys the original topology. Motivated
from the above, this paper uses the FB transform to construct a TBC without
introducing any nonlinear map, and develops an effective BIE solver to compute
the non-quasi-periodic total wavefield. 

To clarify the underlying methodology, we focus only on a periodic array of
sound-soft obstacles with a local perturbation; periodic
structures with different physical properties or different topologies may induce
guided modes and will be studied separately in a future work. First, we develop
a high-accuracy BIE method, following \cite{lu2012high}, to compute the
quasi-periodic Green's function of the background periodic structure for any
quasi-periodic parameter. Integrating these quasi-periodic Green's functions
based on a special quadrature for the inverse FB (IFB) transform, we develop a
fast algorithm to compute the background Green's function and its derivatives
excited by any source point in the periodic structure. By sampling source points
on a sufficiently regular curve surrounding the local perturbation, the exciting
background Green's functions can be exploited to construct a BIE, in terms of
single- and double-layer integral operators on the regular curve. The BIE serves as a
TBC to truncate the unbounded domain. We further derive BIEs on boundaries
of all homogeneities in the truncated domain. The classical kernel-splitting
technique in \cite{colkre13} with spectrally accurate quadratures therein is
readily adapted here to treat the singularities of the background Green's
functions, discretize all BIEs, and establish the final linear system for
the unknown wavefield in the truncated region. By horizontally propagating the quasi-periodic Green's functions, the
background Green's functions in other unit cells, away from the perturbation, can
be easily evaluated by the IFB transform. This allows
us to design two effective procedures, leap and pullback, to compute the
wavefield everywhere in the whole structure. We carry out a number of numerical
experiments to demonstrate that the proposed method can achieve uniform accuracy
no matter the input frequency is half-integer or not, thus outperforming the PML
truncation method \cite{yu2022pml}. Moreover, as the quasi-periodic Green's
functions are evaluated at the boundaries of homogeneities in a single unit cell
for sources on the artificial curve only, the new solver for the
non-quasi-periodic problem exhibits a fully comparable time complexity to that
of one single quasi-periodic problem.

The rest of this paper is organized as follows. In section 2, we present a
mathematical formulation of the scattering problem. In section 3, based on the
Floquet-Bloch theory,  we present a BIE formulation to compute quasi-periodic
Green's functions, and design effective algorithms to compute the background
Green's function. In section 4, we use the background Green’s function to
construct a TBC, that truncates the original problem, and design effective
algorithms to compute the total wavefield. Numerical experiments are carried out
in section 5, followed by some concluding remarks  in section 6 finally.

\section{Problem formulation}
\begin{figure}[htb!]
  \centering
      (a)\includegraphics[width=0.8\textwidth]{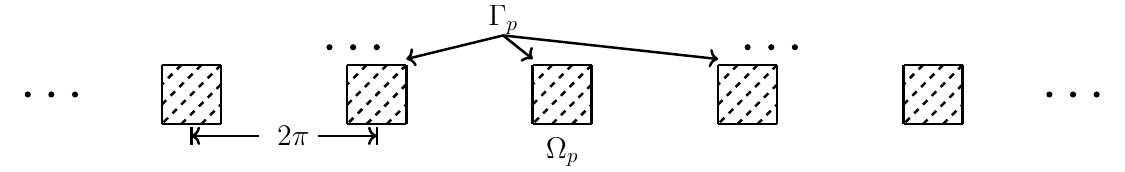}
      (b)\includegraphics[width=0.8\textwidth]{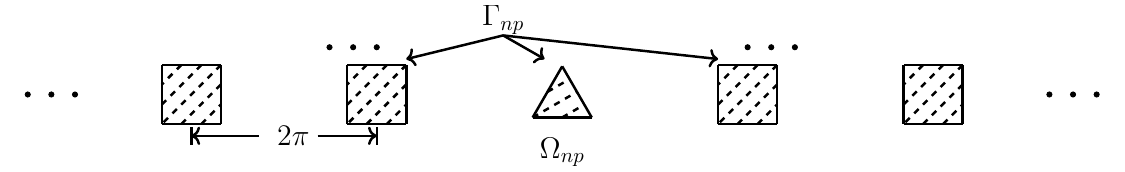}  
   \caption{(a): A $2\pi$-periodic structure with a periodic domain $\Omega_p$; (b): A non-periodic structure with domain $\Omega_{np}$ that locally perturbs $\Omega_p$. }
  \label{fig:1}
\end{figure}

As illustrated in Figure~\ref{fig:1}(a),
in the two dimensional space $\mathbb{R}^2$ with the standard Cartesian
coordinates  $(x_1,x_2)^{T}$, let a periodic array of sound-soft obstacles, with
period $2\pi$ along the $x_1$-direction, be embedded in a homogeneous medium
that occupies the unbounded periodic domain $\Omega_p$. The boundary of
$\Omega_p$, denoted by $\Gamma_p$, consists of the boundaries of all obstacles.
Suppose the periodic structure is locally perturbed, as shown in
Figure~\ref{fig:1}(b), where one obstacle is deformed; 
perturbations with locally different physical properties can be studied 
similarly. We denote the exterior non-periodic domain by $\Omega_{np}$, which is a local
perturbation of $\Omega_p$, and let $\Gamma_{np}$ be its boundary. To
simplify the presentation, we assume that the perturbed region lies in
$|x_1|<\pi$, and all obstacles are separated by the vertical lines $x_1=m\pi,
m\in\mathbb{Z}$.

Suppose an acoustic wave of wavenumber $k>0$ is incident in $\Omega_{np}$ upon the sound-soft
obstacles. The total wavefield $u^{\text{tot}}$ is governed by the 2D Helmholtz
equation: \begin{align}
  \label{eq:1}
  \Delta u^{\text{tot}} + k^2 u^{\text{tot}} &= 0,\quad {\rm in}\ \Omega_{np},
  \end{align}
subject to the homogeneous Dirichlet condition: 
  \begin{align}
  \label{eq:bc}
  u^{\text{tot}} &= 0,\quad {\rm on}\ \Gamma_{np},
\end{align}
where $\Delta=\partial_{x_1}^2+\partial_{x_2}^2\in\Omega_{np}$ is the 2D Laplace operator.  Let $x=(x_1,x_2)^{T}$ be a generic point in $\Omega_{np}$. We will consider two different types of incident waves: (i) a plane incident wave $u^{\text{inc}}(x)=e^{\bi k(\cos\theta x_1-\sin\theta x_2)}$ with the incident angle $\theta\in(0,\pi)$; 
(ii) a cylindrical incident wave $u^{\text{inc}}(x;x^*)=\Phi(x,x^*):=\frac{\bi}{4}H^{(1)}_{0}(k|x-x^*|)$ excited by a source at a point $x^*=(x_1^*,x_2^*)\in \Omega_{np}$.

For the plane-wave incidence (i), let $u^{\text{tot}}_{\text{ref}}(x)$ denote the quasi-periodic total wavefield for the unperturbed structure $\Omega_p$. 
If no propagation modes exist for the periodic structure under consideration, Hu and Kirsch \cite{hukir24} proved that
the perturbed wavefield $u^{\rm og}:= u^{{\rm tot}}-u^{{\rm tot}}_{{\rm ref}}$ satisfies the following Sommerfeld radiation condition (SRC)
\begin{align}
  \label{con:src}
  \lim_{j\to\infty} \left|\left|\frac{\partial u^{\rm og}}{\partial\nu} - \bi k u^{\rm og}\right|\right|_{H^{-1/2}(B_{a_j})} = 0,
\end{align}
where the strictly increasing sequence $\{a_j\}$ is chosen such that
$B_{a_j}:=\{x: x_1^2 + x_2^2=a_j^2\}\subset\Omega_{np}$ and $a_j\to\infty$ as
$j\to\infty$. Note that $u^{\rm tot}_{\rm
ref}$ is easy to compute numerically because of its quasi-periodicity
\cite{lu2012high}.  For the cylindrical incidence (ii), equation \eqref{eq:1} should be modified to 
\begin{align}
  \label{eq:2}
  \Delta u^{\text{tot}} + k^2 u^{\text{tot}} &= -\delta(x-x^*),\quad {\rm in}\quad \Omega_{np}, 
\end{align}
In fact, $u^{\text{tot}}(x;x^*)$ represents the Green's function excited by the
source point $x^*$ for the perturbed structure. Unlike the
plane-wave incidence, $u^{\rm tot}$ satisfies the SRC condition \eqref{con:src}
under the assumption that no propagation modes exist \cite{hukir24}. 

The SRC condition \eqref{con:src} enables us to derive a BIE-based TBC for
$u^{\rm og}$ to truncate the unbounded domain (c.f. Eq.~\eqref{eq:tbcbie}),
using the background Green's function for the unperturbed structure. 

\section{Background Green's function} 
\begin{figure}[ht!]
  \centering
  \includegraphics[width=0.9\textwidth]{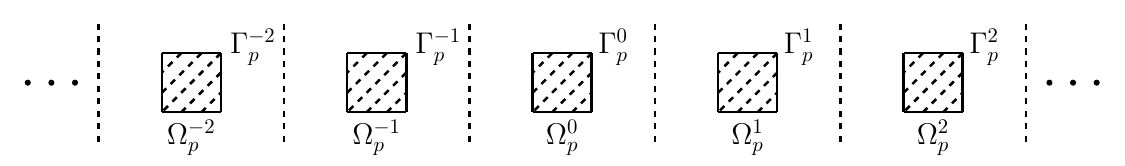}
  \caption{The unbounded homogeneous domain $\Omega_{p}^j$ and the sound-soft obstacle boundary $\Gamma_p^j$ in the $j$-th unit cell for $j\in\mathbb{Z}$.}
  \label{fig:2}
\end{figure}

For a generic source point $x^s\in\Omega_p$, the background Green's function $G(x;x^s)$ is defined as the solution of
\begin{subequations}
\label{eq:_0}
\begin{align}
  \label{eq:1_0}
  \Delta G(x;x^s) + k^2 G(x;x^s) &= -\delta(x-x^s),\quad {\rm in}\ \Omega_p,\\
  \label{eq:2_0}
  G(x;x^s) &= 0,\quad {\rm on}\ \Gamma_p,
\end{align}
\end{subequations}
under the SRC condition \eqref{con:src}. 
As shown in Figure~\ref{fig:2}, let $C_j=((2j-1)\pi, (2j+1)\pi)\times
\mathbb{R}$ and $\Omega_{p}^j := C_j\cap\Omega_p$ be the $j$-th unit cell for
$j\in\mathbb{Z}$. Moreover, we denote the boundary of the $j$-th obstacle by
$\Gamma_{p}^j=C_j\cap\Gamma_p$. 
Without loss of generality, we assume $x^s\in \Omega_{p}^0$, which corresponds
to the perturbed region.

To evaluate the non-quasi-periodic $G(x;x^s)$, we apply the FB transform to express $G(x;x^s)$ as an integral of quasi-periodic Green's
functions.

\subsection{Floquet-Bloch transform and quasi-periodic Green's function}
We recall some basic concepts of the FB transform first. For any
compactly supported function $\varphi\in C^{\infty}_0(\Omega_p)$ with
$\varphi|_{\Gamma_p}=0$, a one-dimensional
FB transform ${\cal J}$ in the $x_1$-direction is defined as
\begin{align}
  \label{fb:1}
  (\mathcal{J}\varphi)(\alpha,x):=\sum_{j\in\mathbb{Z}}\varphi\left(x+\begin{pmatrix}2\pi j\\0\end{pmatrix}\right)e^{2\bi\pi j\alpha},\quad x\in\Omega_{p};\alpha\in\mathbb{R},
\end{align}
where $\alpha$ is called the FB parameter. Clearly, $\psi:={\cal J}\varphi$ belongs to $C^{\infty}([-1/2,1/2]\times\Omega_0)$ and is supported in the $x_2$-direction. 
For a fixed $\alpha\in[-1/2,1/2]$, $\psi$ is $\alpha$-quasiperiodic in $x_1$ with period $2\pi$ in the sense that 

\begin{equation}
\label{eq:qp}  
\left.\psi\left(\alpha,\begin{pmatrix}x_1+2\pi\\x_2\end{pmatrix}\right.\right)=e^{2\bi\pi\alpha}\left.\psi\left(\alpha,\begin{pmatrix}x_1\\x_2\end{pmatrix}\right.\right),\quad \forall x\in\Omega_p,
\end{equation}
and for a fixed $x\in\Omega_p$, $\psi$ is $1$-periodic in $\alpha$, i.e., $\psi(\alpha+1,x)=\psi(\alpha,x)$ for any $\alpha\in[-1/2,1/2]$.
Furthermore, integrating $\psi$ w.r.t. the FB parameter $\alpha$ over $(-1/2,1/2)$ reproduces $\varphi$, inducing the IFB transform ${\cal J}^{-1}$ defined as
\begin{align*}
  (\mathcal{J}^{-1}\psi)(x):=\varphi(x)=\int_{-1/2}^{1/2}\psi(\alpha,x)d\alpha,\quad x\in\Omega_{p}.
\end{align*}
For properties of ${\cal J}$ and ${\cal J}^{-1}$ in weaker functional spaces, we refer readers to \cite{lechleiter2017floquet} for further details. 

Now, let $G^{qp}(\alpha, x;x^s):=({\cal J} G)(\alpha, x;x^s)$, and we shall also denote it by $G^{qp}_\alpha(x;x^s)$ when emphasizing its dependence on $x$. Applying  the FB transform to both sides of \eqref{eq:_0}, $G^{qp}$ is governed by
\begin{subequations}
\label{eqfb}
\begin{align}
  \label{eqfb:1}
  \Delta G^{qp}_\alpha(x;x^s) + k^2 G^{qp}_\alpha(x;x^s) &= -\delta(x-x^s),\quad {\rm in}\ \Omega_{p}^0,\\
  \label{eqfb:2}
  G^{qp}_\alpha(x;x^s) &= 0,\quad {\rm on}\ \Gamma_{p}^0,\\
  \label{eqfb:3}
  G^{qp}_\alpha((\pi,x_2)^T;x^s) &= e^{2\bi\pi\alpha}G^{qp}_\alpha((-\pi,x_2)^{T};x^s),\quad x_2\in\mathbb{R},\\
  \label{eqfb:4}
  \partial_{x_1}G^{qp}_\alpha((\pi,x_2)^T;x^s) &= e^{2\bi\pi\alpha}\partial_{x_1}G^{qp}_\alpha((-\pi,x_2)^T;x^s), \quad x_2\in\mathbb{R},
\end{align}
\end{subequations}

under the condition that $G^{qp}_\alpha$ propagates outwards at infinity. We
remark that $G^{qp}_\alpha(x;x^s)$ is the quasi-periodic Green's function for
the unperturbed background space.  Compared to Eqs.~\eqref{eq:_0},
Eqs.~\eqref{eqfb} are much easier to solve numerically as now we only need to
truncate the $x_2$-direction.  Once $G_\alpha^{qp}$ is obtained, we can evaluate
$G$ by directly applying the IFB transform to $G^{qp}_\alpha$,
and get, 
\begin{equation}
\label{eq:ifbG}
    G(x;x^s) = \int_{-1/2}^{1/2} G^{qp}(\alpha,x;x^s)d\alpha.
\end{equation}

To accurately discretize the integral \eqref{eq:ifbG}, we need to figure out the
regularity of the integrand $G^{qp}$ w.r.t. the FB parameter $\alpha$. According
to \cite{lechleiter2017floquet}, the integrand $G^{qp}$ has two types of
singularities as $\alpha$ ranges over $(-1/2,1/2]$. First, $G^{qp}$ has a
square-root singularity of the form $\sqrt{\kappa\pm\alpha}$ at $\alpha=\mp
\kappa$, where $\kappa\in(-1/2,1/2]$ is chosen such that
$k-\kappa\in\mathbb{N}$. Throughout this paper, we shall use the negative real
axis as the branch cut of the square-root function $\sqrt{\cdot}$ to ensure that
its argument always lies in $[0,\pi)$. The singular integrand can be made
smoother under a rescaling transformation \cite{zhang2018high}, or even made
analytic by simple changes of variables \cite{zhang2022fast}. Second, $G^{qp}$
may have a pole at $\alpha=\alpha_0$ for some $\alpha_0\in(-1/2,1/2]$, which
corresponds to a guided mode of the problem \eqref{eqfb}. In this case, the real
path can be deformed to a new path in the complex plane to bypass $\alpha_0$,
ensuring the analyticity of the integrand \cite{zhang2021numerical}. Then,
high-accuracy quadratures  can be adopted to discretize the tailored integral. 

For the scattering problem under consideration, it is seen from a numerical
perspective that no guided modes exist, though theoretical justifications remain
open.  It suffices to treat the square-root singularity. Suppose
$\kappa\notin\{0,1/2\}$ first. In \cite{zhang2022fast}, the integration domain
$(-1/2,1/2)$ is divided into four sub-intervals, with an appropriate change of
variable applied to each to remove the singularity. In fact, it is enough to
divide it into two sub-intervals. Specifically,
\begin{align}
  \label{ifb:2}
  G=\int_{-1/2}^{1/2}G^{qp}(\alpha)d\alpha=\int_{-\kappa}^{1-\kappa}G^{qp}(\alpha)d\alpha=(\int_{-\kappa}^{\kappa}+\int_{\kappa}^{1-\kappa})G^{qp}(\alpha)d\alpha,
\end{align}
where the arguments $x$ and $x^s$ are omitted for brevity.
For the first integral on the r.h.s, to eliminate the square-root singularities
$\sqrt{\kappa\pm \alpha}$ at $\alpha=\mp\kappa$, we introduce the substitution
$\alpha=\kappa\cos\theta$ and obtain
\begin{align}
\label{ifb:2_1}
  \int_{-\kappa}^{\kappa}G^{qp}(\alpha)d\alpha=\kappa\int_{0}^{\pi}G^{qp}(\kappa\cos\theta)\sin\theta d\theta. 
\end{align}
For the second integral, we use $\alpha=\frac{(2\kappa-1)\cos\theta+1}{2}$ and get
\begin{align}
\label{ifb:2_2}
  \int_{\kappa}^{1-\kappa}G^{qp}(\alpha)d\alpha=\frac{1-2\kappa}{2}\int_{0}^{\pi}G^{qp}\left(\frac{(2\kappa-1)\cos\theta+1}{2}\right)\sin\theta d\theta.
\end{align}
For an even integer $n$, we adopt the $n/2$-point Gauss-Legendre quadrature to discretize each integral in \eqref{ifb:2_1} and \eqref{ifb:2_2}, and get 
\begin{equation}
\label{eq:disG}
    G \approx \sum_{i=1}^{n} G^{qp}(\alpha_i)w_i=\sum_{i=1}^{n} G^{qp}_{\alpha_i} w_i, 
\end{equation}
where $w_i$ and $\alpha_i$ are the weights and the $n$ distinct nodes in $(-\kappa,1-\kappa)$, respectively. We note that $n$ does not have to be proportional to the wavenumber $k$, as $\alpha$ is confined to $(-1/2,1/2]$.  For $\kappa\in\{0,1/2\}$, we do not need to split the integration domain, and the corresponding change of variable is similar. 

\subsection{A BIE formulation for Problem \eqref{eqfb}}

In this subsection, we present a BIE formulation for computing $G^{qp}_\alpha$ for $\alpha=\alpha_i$ as in \eqref{eq:disG}. It follows a previous work \cite{lu2012high} that computes quasi-periodic wavefields for incident plane waves. We note that \eqref{eq:disG} suggests that $\alpha$ can be close to $\pm \kappa$ or $1-\kappa$, corresponding to the Rayleigh anomaly, making it challenging to use quasi-periodic fundamental solutions to compute $G^{qp}_\alpha$ \cite{bruno20}. 
The key step of \cite{lu2012high} is to formulate an effective boundary condition that truncates the unbounded domain $\Omega_{p}^0$.
Write $\Omega_{p,H}^{0}:=\Omega_{p}^0\cap\{x:|x_2|\leq H\}$ for $H>|x_2^s|$.  
We shall frequently omit the source point $x^s$ to write
$G^{qp}_\alpha(x)$ for brevity. 
\begin{figure}[ht!]
\centering
\includegraphics[width=0.5\textwidth]{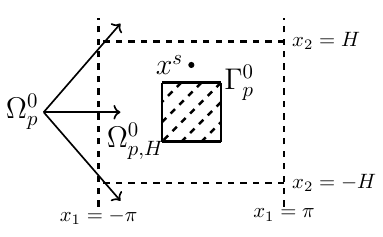}  
\caption{The periodic structure in the $0$-th unit cell $\Omega_p^0$: $x^s$ represents a generic source point, TBCs \eqref{top:con} and \eqref{bottom:con} are established on $x_2=\pm H$, and $\Omega_{p,H}^0$ denotes the truncated domain. }
\label{fig:3.1}
\end{figure}

Owing to its quasi-periodicity, $G^{qp}_\alpha$ admits the following Rayleigh expansions 
\begin{align}
\label{raylei}
  G^{qp}_\alpha(x)=\sum\limits_{l\in \mathbb{Z}}B_{l}^\pm e^{\bi(\alpha_lx_1\pm \beta_l(x_2-H))}, \pm x_2\geq H,
\end{align}
where $\alpha_l = \alpha + l$, $\beta_l = \sqrt{k_0^2-\alpha_l^2}$, and
\begin{align}
\label{fourierco}
    B_l^\pm=\frac{1}{2\pi}\int_{-\pi}^{\pi}G^{qp}_\alpha((x_1,\pm H)^{\rm T})e^{-\bi\alpha_l \tx_1}d\tx_1.
\end{align} 
On $x_2=H$, we get 
\begin{align}
  \label{top:con}
  \partial_{x_2}G^{qp}_\alpha\left( \begin{pmatrix}x_1\\H\end{pmatrix}\right)={\cal B}_\alpha G^{qp}_\alpha\left(\begin{pmatrix}x_1\\H\end{pmatrix}\right),
\end{align}
where the Dirichlet-to-Neumann operator ${\cal B}_\alpha$ is defined as: for any $\alpha$-quasiperiodic function $f\in C^{\infty}([-\pi,\pi])$, 
\begin{align}
  \label{tbc}
({\cal B}_\alpha f)(x)=\frac{\bi}{2\pi}\sum\limits_{l=-\infty}^{\infty}\beta_l\left(\int_{-\pi}^{\pi}f(\tilde{x})e^{-\bi\alpha_l\tilde{x}_1}d\tilde{x}_1\right) e^{\bi\alpha_l x_1}.
\end{align}
It can be verified that ${\cal B}_\alpha
e^{\bi\alpha_nx_1}=\bi\beta_ne^{\bi\alpha_nx_1}$ for all $n$. Similarly, on
$x_2=-H$, 
\begin{align}
  \label{bottom:con}
  \partial_{x_2}G^{qp}_\alpha\left(\begin{pmatrix}x_1\\-H\end{pmatrix}\right)=-{\cal B}_\alpha G^{qp}_\alpha\left(\begin{pmatrix}x_1\\-H\end{pmatrix}\right).
\end{align}
Equations (\ref{top:con}) and (\ref{bottom:con}) serve as TBCs for the unknown wavefield $G^{qp}_\alpha$ on $x_2=\pm H$.

The scattered wave $G^{sc}_\alpha(x;x^s):={G^{qp}_\alpha}(x;x^s)-\Phi(x;x^s)$ satisfies the homogeneous Helmholtz equation 
\begin{align}
  \label{eq:11}
  \Delta G^{sc}_\alpha + k^2 G^{sc}_\alpha &= 0,\quad {\rm in}\ \Omega_{p,H}^{0}.
\end{align} 
Thus, $G^{sc}$ admits the following Green's representation formula
\begin{align}
  \label{rp:1}
G^{sc}_\alpha(x)=\int_{\partial\Omega_{p,H}^0}\left[\Phi(x;y)\partial_{\nu}G^{sc}_\alpha(y)-\partial_{\nu}\Phi(x;y)G^{sc}_\alpha(y)\right]ds(y),\quad x\in\Omega_{p,H}^{0},
\end{align} 
where $\nu$ denotes the outer unit normal to the boundary $\partial\Omega_{p,H}^0$. 
As $x$ approaches $\partial\Omega_{p,H}^0$, one obtains the following boundary integral equation \cite{lulu14}
\begin{align}
  \label{bie:1}
  (\mathcal{K}-\mathcal{K}_0[1])[G^{sc}_\alpha](x)=\mathcal{S}[\partial_{\nu}G^{sc}_\alpha](x),
\end{align} 
for $x\in \partial\Omega_{p,H}^0$. Here, the boundary integral operators are defined as
\begin{align}
  \label{bie:o1}
  \mathcal{S}[\phi](x)&=2\int_{\partial\Omega_{p,H}^0}\Phi(x;y)\phi(y)ds(y),\\
  \label{bie:o2}
  \mathcal{K}[\phi](x)&=2\text{P.V.}\int_{\partial\Omega_{p,H}^0}\partial_{\nu}\Phi(x;y)\phi(y)ds(y),\\
  \label{bie:o3}
  \mathcal{K}_0[\phi](x)&=2\text{P.V.}\int_{\partial\Omega_{p,H}^0}\partial_{\nu}\Phi_0(x;y)\phi(y)ds(y),
\end{align}
where $\Phi_0(x;y)=\frac{1}{2\pi}\log|x-y|$ is the fundamental solution of the two-dimensional Laplacian,  and $\text{P.V.}\int$ denotes the Cauchy principal integral. 
This induces a Neumann-to-Dirichlet (NtD) operator 
\begin{equation}
\label{eq:NtD}
    \mathcal{N}=(\mathcal{K}-\mathcal{K}_0[1])^{-1}\mathcal{S},
\end{equation}
which maps $\partial_{\nu}G^{sc}_\alpha$ to $G^{sc}_\alpha$ on the boundary $\partial\Omega_{p,H}^0$.

\begin{figure}[htb]
\centering
\includegraphics[width=0.5\textwidth]{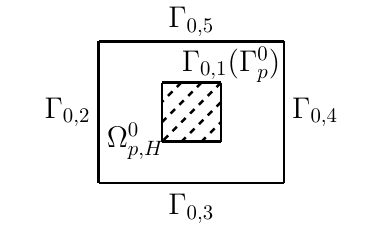}
\caption{The inner, left, bottom, right and top parts $\Gamma_{0,j}$ for $j=1,\cdots,5$.}
\label{fig:4.1}
\end{figure}

We now formulate the problem in the form of abstract operators within the bounded rectangular domain $\Omega_{p,H}^0$. 
As shown in Figure~\ref{fig:4.1}, we denote the inner, left, bottom, right and top parts of
$\partial\Omega_{p,H}^0$ as $\Gamma_{0,i}$ for $1\leq i\leq 5$, respectively;
note that $\Gamma_{0,1} = \Gamma_p^0$. 
Write
$\psi_i=G^{sc}_\alpha|_{\Gamma_{0,i}}$,
$\partial_{\nu}\psi_i=\partial_{\nu}G^{sc}_\alpha|_{\Gamma_{0,i}}$,
$\Phi_i=\Phi|_{\Gamma_{0,i}}$, and
$\partial_{\nu}\Phi_i=\partial_{\nu}\Phi|_{\Gamma_{0,i}}$. Thus, 
\begin{align}
  \label{ab:2}
  \left.\left[\begin{array}{c}{{\psi}_{1}}\\{{\psi}_{2}}\\{{\psi}_{3}}\\{{\psi}_{4}}\\{\psi}_{5}\\\end{array}\right.\right]=
  \mathcal{N}\left[\begin{array}{c}{{\partial_{\nu}\psi}_{1}}\\{{\partial_{\nu}\psi}_{2}}\\{{\partial_{\nu}\psi}_{3}}\\{{\partial_{\nu}\psi}_{4}}\\{\partial_{\nu}\psi}_{5}\\\end{array}\right],
\end{align}
On the other hand, the quasi-periodic conditions \eqref{eqfb:3} and \eqref{eqfb:4}, the top and bottom conditions \eqref{top:con} and \eqref{bottom:con}, and the Dirichlet condition \eqref{eqfb:2} lead to a matrix system 
\begin{align}
  \label{ab:1}
  {\cal P}_\alpha\left(\left.\left[\begin{array}{c}{{\psi}_{1}}\\{{\psi}_{2}}\\{{\psi}_{3}}\\{{\psi}_{4}}\\{\psi}_{5}\\
  \end{array}\right.\right]+\left.\left[\begin{array}{c}{{\Phi}_{1}}\\{{\Phi}_{2}}\\{{\Phi}_{3}}\\{{\Phi}_{4}}\\{\Phi}_{5}\\
  \end{array}\right.\right]\right)=
  {\cal Q}_\alpha\left(\left[\begin{array}{c}{{\partial_{\nu}\psi}_{1}}\\{{\partial_{\nu}\psi}_{2}}\\{{\partial_{\nu}\psi}_{3}}\\{{\partial_{\nu}\psi}_{4}}\\{\partial_{\nu}\psi}_{5}\\
  \end{array}\right]+\left.\left[\begin{array}{c}{{\partial_{\nu}\Phi}_{1}}\\{{\partial_{\nu}\Phi}_{2}}\\{{\partial_{\nu}\Phi}_{3}}\\{{\partial_{\nu}\Phi}_{4}}\\{\partial_{\nu}\Phi}_{5}\\
  \end{array}\right.\right]\right),
\end{align}
where  ${\cal I}$ denotes the abstract identity operator, and
\begin{align}
  \label{ab:3}
  {\cal P}_\alpha=\left.\left[\begin{matrix}
  0&\mathcal{I}&0&-e^{2\bi\pi\alpha}\mathcal{I}&0\\
  0&0&0&0&0\\
  0&0&-{\cal B}_\alpha &0&0\\
  0&0&0&0&{\cal B}_\alpha\\
  \mathcal{I}&0&0&0&0\\
  \end{matrix}\right.\right],\quad
  {\cal Q}_\alpha=\left.\left[\begin{matrix}
  0&0&0&0&0\\
  0&\mathcal{I}&0&e^{2\bi\pi\alpha}\mathcal{I}&0\\
  0&0&{\cal I}&0&0\\
  0&0&0&0&{\cal I}\\
  0&0&0&0&0\\
  \end{matrix}\right.\right], 
\end{align}
Combining \eqref{ab:1} and \eqref{ab:2}, 
\begin{align}
\label{ab:4}
    ({\cal P}_\alpha\mathcal{N}-{\cal Q}_\alpha)\left[\begin{array}{c}{{\partial_{\nu}\psi}_{1}}\\{{\partial_{\nu}\psi}_{2}}\\{{\partial_{\nu}\psi}_{3}}\\{{\partial_{\nu}\psi}_{4}}\\{\partial_{\nu}\psi}_{5}\\
  \end{array}\right]=
    {\cal Q}_\alpha\left.\left[\begin{array}{c}{{\partial_{\nu}\Phi}_{1}}\\{{\partial_{\nu}\Phi}_{2}}\\{{\partial_{\nu}\Phi}_{3}}\\{{\partial_{\nu}\Phi}_{4}}\\{\partial_{\nu}\Phi}_{5}\\
  \end{array}\right.\right]-{\cal P}_\alpha\left.\left[\begin{array}{c}{{\Phi}_{1}}\\{{\Phi}_{2}}\\{{\Phi}_{3}}\\{{\Phi}_{4}}\\{\Phi}_{5}\\
  \end{array}\right.\right].
\end{align}

On solving \eqref{ab:4}, we obtain $\partial_\nu G^{sc}_\alpha$ on
$\partial\Omega_{p,H}^0$. Then, $G^{sc}_\alpha={\cal N} \partial_\nu
G^{sc}_\alpha$ on $\partial\Omega_{p,H}^0$. By Green's representation formula
\eqref{rp:1}, $G^{sc}_\alpha$ and hence the total field $G^{qp}_\alpha =
G^{sc}_\alpha + \Phi$ become available in $\Omega_{p,H}^{0}$. In the exterior region
$\Omega_{p}^0\backslash\overline{\Omega_{p,H}^0}$, we apply \eqref{fourierco} to
obtain the Rayleigh coefficients $B_l^\pm$ first and then use \eqref{raylei} to
evaluate $G^{qp}_\alpha$. The following subsection is devoted to accurate approximations of ${\cal N}$ and ${\cal B}_\alpha$. 

\subsection{Numerical approximations of ${\cal N}$ and ${\cal B}_\alpha$}

Following \cite{lulu14}, we apply the Nystr\"om method to discretize the three
integral operators in \eqref{bie:o1}-\eqref{bie:o3}. For completeness and
for the later treatment of the TBC proposed in
section 4.1, we present its brief idea below. We study the single-layer
operator ${\cal S}$ only; the other two operators are treated similarly.

As indicated in Figure~\ref{fig:4.1}, 
$\partial\Omega_{p,H}^0$ consists of two well-separated closed curves, a rectangle $\gamma_1:=\partial([-\pi,\pi]\times[-H,H])=\bigcup\limits_{j=2}^5\Gamma_{0,j}$ and the obstacle boundary $\gamma_2:=\Gamma_{0,1}$. To simplify the presentation, we assume $x\in\gamma_1$ for the moment, and write
\begin{align}
  \label{S:1}
  \mathcal{S}[\partial_\nu G_\alpha^{sc}](x)&=\int_{\gamma_1}\Phi(x;y)\partial_\nu G_\alpha^{sc}(y) ds(y)+\int_{\gamma_2}\Phi(x;y)\partial_\nu G_\alpha^{sc}(y) ds(y)\nonumber\\
  &:={\cal S}^{(1)}[\partial_\nu G_\alpha^{sc}](x) + {\cal S}^{(2)}[\partial_\nu G_\alpha^{sc}](x).
\end{align} 
We consider the weakly singular operator ${\cal S}^{(1)}$ only.

Suppose $\gamma_1$ is parameterized by $x(s)=\{(x_1(s),x_2(s))|0\leq s\leq L\}$ for the arc-length parameter $s$. Thus, we can parameterize ${\cal S}^{(1)}$ as 
\begin{equation}
\label{eq:S11}
    {\cal S}^{(1)}[\partial_\nu G_\alpha^{sc}](x(s)) = \int_{0}^L \Phi(x(s), x(s')) \partial_\nu G_\alpha^{sc}(x(s'))ds'.
\end{equation}
We introduce a scaling function $s=\eta(t), t\in[0,2\pi]$ to smoothen the
discontinuous density function $\partial_\nu G_\alpha^{sc}(x(s))$. Suppose
$s_b=\eta(t_b)$ and $s_e=\eta(t_e)$ correspond to two consecutive corners with
$0\leq t_b<t_e\leq 2\pi$ and $0\leq s_b<s_e\leq L$. We define \cite[Eq.~(3.104)]{colkre13}
%
\begin{align}
  \label{cov:1}
  \eta(t)=\frac{s_b\eta_1^p+s_e\eta_2^p}{\eta_1^p+\eta_2^p},\quad t\in[t_b,t_e],
\end{align} 
where $p$ is a positive integer, and
\begin{align}
  \label{cov:2}
  \eta_1=(\frac{1}{2}-\frac{1}{p})\xi^3+\frac{\xi}{p}+\frac{1}{2},\quad \eta_2=1-\eta_1, \quad \xi=\frac{2t-(t_b+t_e)}{t_e-t_b}. 
\end{align} 
Using a uniform mesh of $N$ points $\{t_j=2\pi
j/N\}_{j=0}^{N-1}$ on $[0,2\pi]$ for an even integer $N>0$,
$s_j=\eta(t_j)$ generates a graded mesh with grid points clustering near corners
and uniformly spaced away from them; in our implementation, corners are chosen
as part of the grid points and the mesh points on $\Gamma_{0,4}$ should be the
translate of those on $\Gamma_{0,2}$. To simplify the notation, we use $x(t)$
and $x'(t)$ to denote $x(\eta(t))$ and $x'(\eta(t))\eta'(t)$, respectively.




Now, equation \eqref{eq:S11} becomes
\begin{align}
  \label{SS:1}
  {\cal S}^{(1)}[\partial_\nu G_\alpha^{sc}](x(t))=\int_{0}^{2\pi}S(t,\tau)\phi_\alpha^{sc}(\tau)d\tau, 
\end{align}
where $\phi_\alpha^{sc}(\tau) = \partial_\nu G_\alpha^{sc}(\tau)\eta'(\tau)$ and 
$S(t,\tau):= \frac{\bi}{2}H_0^{(1)}(k|x(t)-x(\tau)|)$; here we use the
scaled normal derivative $\phi_\alpha^{sc}$ in place of $\partial_\nu
G_\alpha^{sc}$ because it can significantly reduce the condition number of the
approximate matrix of ${\cal S}^{(1)}$ given below \cite{lulu12}. The kernel function $S$ can be split as \cite{colkre13} 
\begin{align}
  S(t,\tau)= S_1(t,\tau)\ln(4\sin^2\frac{t-\tau}{2})+S_2(t,\tau),
\end{align}
where 
\begin{align}
  \label{ks:2}
  S_1(t,\tau)&= -\frac{1}{2\pi}J_0(k|x(t)-x(\tau)|),\\
  \label{ks:3}
  S_2(t,\tau)&= \begin{cases}
      S(t,\tau)-S_1(t,\tau)\ln(4\sin^2\frac{t-\tau}{2}),& t\neq \tau,\\
      \displaystyle
      \frac{\bi}{2}-\frac{C}{\pi}-\frac{1}{\pi}\ln(\frac{k}{2}|x'(t)|),& t=\tau,
  \end{cases}
\end{align}
and $C$ denotes Euler's constant. Then, we approximate
\begin{align}
  \label{ks:6}
{\cal S}^{(1)}[\partial_\nu G_\alpha^{sc}](x(t_j))\approx\sum\limits_{k=0}^{N-1}[R_{|k-j|}^{(N)}S_1(t_j,t_k)+\frac{2\pi}{N}S_2(t_j,t_k)]\phi_\alpha^s(t_{k}), j=0,\cdots, N-1,
\end{align}
where 
\[
R_n^{(N)}:=-\frac{4\pi}{N}\sum\limits_{m=1}^{N/2}\frac{1}{m}\cos \frac{2m n \pi}{N}-\frac{4\pi}{N^2}(-1)^n.
\]
This gives rise to 
\begin{equation}
 {\cal S}^{(1)}[\partial_\nu G_\alpha^{sc}]
 \left[
 \begin{matrix}
     x(t_0)\\
     x(t_1)\\
     \vdots\\
     x(t_{N-1})
 \end{matrix}
 \right]= {\bf S}_{11} \left[
 \begin{matrix}
     \phi_\alpha^s(t_0)\\
     \phi_\alpha^s(t_1)\\
     \vdots\\
     \phi_\alpha^s(t_{N-1})
 \end{matrix}
 \right],
\end{equation}
for an $N\times N$ matrix $\textbf{S}_{11}$, where the l.h.s. represents a column vector of elements ${\cal S}^{(1)}[\partial_\nu G_\alpha^{sc}](x(t_j))$.

To discretize ${\cal S}^{(2)}[\partial_\nu G_\alpha^{sc}](x(t_j))$, we can use the same approach, in fact simpler due to $x(t_j)\notin\gamma_2$, to obtain an $N\times M$ matrix ${\bf S}_{12}$, where $M$ is the number of grid points on $\gamma_2$. Note that the $M$ grid points are uniformly spaced in the parameter space when $\gamma_2$ is smooth, and form a graded mesh otherwise. Now consider \eqref{S:1} for $x\in\gamma_2$. Using the same grid points on $\Gamma_0$, we can obtain an $M\times N$ matrix $\textbf{S}_{21}$ for the operator ${\cal S}^{(1)}$ and an $M\times M$ matrix $\textbf{S}_{22}$ for the operator ${\cal S}^{(2)}$.  Consequently, the matrix ${\bf S}=[{\bf S}_{ij}]_{i,j=1,2}$ represents an accurate approximation of the single-layer operator ${\cal S}$ on $\gamma_1\cup\gamma_2$. 

Applying a similar procedure to ${\cal K}$ and ${\cal K}_0$ as above, \eqref{bie:1} can be approximated by
\begin{equation}
    (\textbf{K}-{\bf D}){\bm G}_\alpha^{sc} = {\bf S}{\bm \phi}_\alpha^{sc},
\end{equation}
where ${\bf K}$ is an $(M+N)\times (M+N)$ matrix approximating ${\cal K}$, ${\bf D}$ is a diagonal matrix with diagonal elements approximating ${\cal K}_0[1]$, and ${\bm G}_\alpha^{sc}$ (${\bm \phi}_\alpha^{sc}$) represents a column vector of $G_\alpha^{sc}$ ($\phi_\alpha^{sc}$) at the $M+N$ grid points on $\gamma_1\cup\gamma_2$. This gives rise to an $(M+N)\times (M+N)$ matrix 
\begin{equation}
\label{eq:ntd}
    \textbf{N}=(\textbf{K}-{\bf D})^{-1}\textbf{S},
\end{equation}
approximating the NtD operator ${\cal N}$.

Now, suppose the top segment of $\gamma_1$, namely $\Gamma_{0,5}$, is parameterized in $t$ as $x(t)$ for $0\leq t\leq T$, and the graded mesh points on it are $\{t_k\}_{k=0}^{N_0-1}$ for some integer $N_0>0$. For the top boundary condition \eqref{top:con},  
we can apply the trapezoidal rule to the integrals in \eqref{tbc} and get \cite{lu2012high} 
\begin{align*}
  \label{dtb}
 \phi_\alpha^{qp}(x(t_j)):=&\eta'(t_j)\partial_{\nu}G_\alpha^{qp}(x(t_j))  = \eta'(t_j) ({\cal B}_\alpha G_\alpha^{qp})(x(t_j))\\
  \approx&\frac{\bi}{2\pi}\sum\limits_{l=-J_0}^{J_0}\eta'(t_j)\beta_{l}e^{\bi\alpha_l x_1(t_j)}\int_{0}^{T}G_\alpha^{qp}(x(\tau))e^{-\bi\alpha_l x_1(\tau)}\eta'(\tau)d\tau\\
  \approx&\sum\limits_{k=0}^{N_0-1}\left[\frac{\bi T\eta'(t_j)}{2\pi N_0}\sum\limits_{l=-J_0}^{J_0}\beta_{l}e^{\bi\alpha_l x_1(t_j)}e^{-\bi\alpha_l x_1(t_k)}\eta'(t_k)\right]G_\alpha^{qp}(x(t_k)),
\end{align*}
for $0\leq j\leq N_0-1$. Since the Rayleigh coefficients $B_l$ decay
exponentially with $l$, a moderately large integer $J_0$ ensures 
sufficient accuracy in practice. The bracketed coefficients can be
rapidly evaluated via the nonuniform fast Fourier transform, and form an $N_0\times N_0$ matrix
$\textbf{B}_\alpha$ such that 
\begin{equation} 
\label{eq:Balp}
{\bm \phi}_\alpha^{qp} = {\bf B}_\alpha {\bm G}_\alpha^{qp}, 
\end{equation} 
where ${\bm \phi}_\alpha^{qp}$ is a column vector of the $N_0$ elements
$\phi_\alpha^{qp}(x(t_j))$, etc. Roughly speaking, ${\bf B}_\alpha$ can be regarded as
an approximation of the scaled operator $\eta'{\cal B}_\alpha$.

\subsection{Evaluating Green's function and its derivatives}

With the two approximate matrices ${\bf N}$ and ${\bf B}_\alpha$, we are ready to evaluate the Green's function $G(x;x^s)$ and its derivatives in all cells $\Omega^j_{p}$, which are the key ingredients in solving the original scattering problem \eqref{eq:1}-\eqref{eq:2}. Let $\Omega^j_{p,H}:=\Omega^j_p\cap\{x:|x_2|\leq H\}$ and $\partial\Omega^j_{p,H}$ be its boundary for any $j\in\mathbb{Z}$.

We first discuss the Green's function $G$ in $\Omega_p^0$. Equation \eqref{ab:4}
can be approximated by the following linear system \begin{equation}
\label{eq:linsys} {\bm A}(\alpha){\bm \phi}_\alpha^{sc}(x^s) = {\bm
b}(\alpha;x^s), \end{equation} where the $(M+N)\times(M+N)$ matrix ${\bm A}$ is
related to the operator ${\cal P}_\alpha {\cal N}-{\cal Q}_\alpha$ and the
nonzero $(M+N)\times 1$ column vector ${\bm b}$ approximates the r.h.s. of
\eqref{ab:4}. Note that we have multiplied the third and fourth row of
\eqref{ab:4} with $w'$ to  before establishing \eqref{eq:linsys}. Solving the
linear system, we get ${\bm \phi}_\alpha^{sc}$ and then ${\bm G}_\alpha^{sc} =
{\bm N} {\bm \phi}_\alpha^{sc}$ on the boundary $\partial\Omega_{p,H}^0$. 

For any observation point $x\in \Omega_{p,H}^0$, using the trapezoidal rule with the mesh points on the boundary $\partial \Omega_{p,H}^0$ as the nodes to discretize the integrals in the two representation formulae \eqref{rp:1}, we can evaluate $G_\alpha^{sc}(x;x^s)$ from ${\bm \phi}_\alpha^{sc}$ and ${\bm G}_\alpha^{sc}$.  With the IFB transform \eqref{ifb:2} and the quadrature \eqref{eq:disG}, we can evaluate the background Green's function $G(x;x^s)$ and its gradient $\nabla_x G(x;x^s)$ for any $x\in \Omega^0_{p,H}$. The above procedure can be simply understood as numerical approximations of 
\begin{align}
\label{eq:appG}
    G(x;x^s) =& \Phi(x;x^s) \nonumber\\
    &+\int_{-\kappa}^{1-\kappa}d\alpha \int_{\partial\Omega_{0,H}^0} \left[\Phi(x;z)\partial_{\nu(z)}G_\alpha^{sc}(z;x^s) - \partial_{\nu(z)}\Phi(x;z)G_\alpha^{sc}(z;x^s) \right]ds(z),
\end{align}
and its $x$-gradient. Suppose we need to evaluate $G$ at $N_o$ observation points $\{x^o_j\}_{j=1}^{N_o}$ in $\Omega_{p,H}^0$ for the fixed source point $x^s$. As $\partial \Omega_{p,H}^0$ is discretized by $M+N$ points and the interval $(-\kappa, 1-\kappa)$ by $n$ points, it is not hard to find that the computational complexity above is ${\cal O}(N_o n(M+N))$. If we further require $G$ at $N_s$ source points $\{x^s_k\}_{k=1}^{N_s}$, then the complexity becomes ${\cal O}(N_o N_s n(M+N))$, which is too costly. Nevertheless, by interchanging the order of integrations in \eqref{eq:appG}, 
\begin{align}
\label{eq:appG2}
    G(x;x^s) =& \Phi(x;x^s) 
    +\int_{\partial\Omega_{0,H}^0} \Bigg[\Phi(x;z)\int_{-\kappa}^{1-\kappa}\partial_{\nu(z)}G_\alpha^{sc}(z;x^s)d\alpha \nonumber\\
    &- \partial_{\nu(z)}\Phi(x;z)\int_{-\kappa}^{1-\kappa}G_\alpha^{sc}(z;x^s)d\alpha  \Bigg]ds(z)\nonumber\\
    =& \Phi(x;x^s) + 
    \int_{\partial\Omega_{0,H}^0} \left[\Phi(x;z) \partial_{\nu(z)}G^{sc}(z;x^s) - \partial_{\nu(z)}\Phi(x;z) G^{sc}(z;x^s)  \right]ds(z).
\end{align}
Eq.~\eqref{eq:appG2} evaluates $\partial_{\nu(z)}G^{sc}(z;x^s)$ and $G^{sc}(z;x^s)$ for
$z\in\partial\Omega_{p,H}^0$ first, significantly reducing the complexity
to ${\cal O}((N_o+n) N_s (M+N))$. 
This makes the numerical evaluations almost $n$ times faster!

According to \eqref{eqfb} and $\nabla_{x^s} G_\alpha^{sc}(x;x^s)=\nabla_{x^s}
G_\alpha^{qp}(x;x^s)-\nabla_{x^s}\Phi(x;x^s)$, by replacing $G_\alpha^{qp}$ by
$\nabla_{x^s} G_\alpha^{qp}$ and the r.h.s. of \eqref{eqfb:1} by
$-\nabla_{x^s}\delta(x-x^s)$, we can evaluate $\nabla_{x^s} G_\alpha^{qp}$ by
the same method in Section 3.2. As for $x\in \Omega_p^0\backslash
\overline{\Omega_{p,H}^0}$, we use $G_\alpha^{qp}$ on $x_2=\pm H$ to get the
Rayleigh coefficients by \eqref{fourierco} and then \eqref{raylei} to compute
$G_\alpha^{qp}$ and its derivatives. The IFB transform
\eqref{eq:disG} then applies to compute $G$ and its derivatives in
$\Omega_p^0\backslash\overline{\Omega_{p,H}^0}$; we omit the details here. The
above procedure is summarized in Algorithm \ref{alg:evalG} below. 
\begin{algorithm}[htbp]
  \caption{Evaluating the Green's function and its gradient in $\Omega_{p}^0$.} 
  \label{alg:evalG}
  \begin{algorithmic}[1]
    \Require
      $n$ quadrature nodes $\{\alpha_i\}_{i=1}^n\subset(-\kappa,1-\kappa)$;
      $M+N$ mesh points $\{x_l\}_{l=1}^{M+N}\subset\partial\Omega_{p,H}^0$;
      $N_o$ observation points $\{x^o_j\}_{j=1}^{N_o}\subset\Omega_{p,H}^0$;
      $N_s$ source points $\{x_k^s\}_{k=1}^{N_s}\subset\Omega_{p,H}^0$.
       \State Approximate the operator ${\cal N}$ by ${\bf N}$ as in \eqref{eq:ntd}, and the operator ${\cal B}_{\alpha_i}$ by ${\bf B}_{\alpha_i}$ as in \eqref{eq:Balp} for $\alpha=\alpha_i$;
         \For{$i = 1$; $i<n$; $i++$ }
      \State Solve the linear system $\mathbf{A}(\alpha_i){\bm \phi}_{\alpha_i}^{sc}(x_k^s)=\mathbf{b}(\alpha_i;x_k^s)$ according to \eqref{eq:linsys} for ${\bm \phi}_{\alpha_i}^{sc}(x_k^s)$ on $\partial\Omega_{p,H}^0$; 
     \EndFor
      \State Evaluate $G^{sc}(x_l;x_k^s)$, $\partial_\nu G^{sc}(x_l;x_k^s)$ and their gradients w.r.t $x^s$ by \eqref{eq:disG};
      \State Evaluate $G(x^o_j;x_k^s)$ and $\nabla_{x^s} G(x^o_j;x_k^s)$ via \eqref{eq:appG2};
      \State Evaluate $G(x;x_k^s)$ and $\nabla_{x^s} G(x;x_k^s)$ by \eqref{raylei} and then \eqref{eq:disG} for any $x\in \Omega^0_{p}\backslash\overline{\Omega_{p,H}^0}$. 
  \end{algorithmic}
\end{algorithm}

As $\mathbf{A}(\alpha_i)$ depends on $\alpha_i$ only, we can compute the LU
decomposition of $\mathbf{A}(\alpha_i)$ just once for all source points $x_k^s$.
Thus, the computational complexity for solving the $nN_s$ linear systems becomes ${\cal O}((M+N)^2 n (M+N+N_s))$. In
summary, the computational complexity of Algorithm~\ref{alg:evalG} for Steps 1 through 6 is
\[
    {\cal O}((M+N)^3) + {\cal O}((M+N)^2 n (M+N+N_s)) + {\cal O}(N_o N_s n(M+N)),
\]
where the first term accounts for computing ${\bf N}$, the second for the $nN_s$ linear
systems, and the last for evaluating $G$ and its gradient. In practice, the
prefactor of the first term is much larger than the second term, making it the
most costly to approximate the NtD map ${\cal N}$ during the whole procedure. In
other words, the speed of Algorithm~\ref{alg:evalG} for evaluating the
background Green's function $G$ will not be significantly slower than that for
the simpler quasi-periodic Green's function $G^{qp}_\alpha$, as both require
approximating ${\cal N}$!
\begin{algorithm}[H]
  \caption{Evaluating the Green's function and its gradient in $\Omega_{p}^{j}$ for $j\neq 0$.} 
  \label{alg:evalG2}
  \begin{algorithmic}[1]
    \Require
      $M+N$ mesh points $\{x_l\}_{l=1}^{M+N}$ on $\partial\Omega_{p,H}^0$;
      $M+N$ mesh points $\{x_l^{(j)}\}_{l=1}^{M+N}$ on $\partial\Omega_{p,H}^{j}$; 
      $N_s$  source points $\{x_k^s\}_{k=1}^{N_s}\subset\Omega_{p,H}^0$.
       \State Compute ${\bm \phi}_{\alpha_i}^{sc}(x_k^s)$, ${\bm
       G}_{\alpha_i}^{sc}(x_k^s)$ and their gradients w.r.t $x^s$ on $\partial\Omega_{p,H}^{0}$ by Algorithm \ref{alg:evalG}; 
       \State Compute ${\bm \phi}_{\alpha_i}^{sc,(j)}(x_k^s)$, ${\bm
       G}_{\alpha_i}^{sc,(j)}(x_k^s)$ and their gradients w.r.t $x^s$ on $\partial\Omega_{p,H}^{j}$ by \eqref{eq:qp}; 
       \State Evaluate $G^{sc}(x_l^{(j)};x_k^s)$, $\partial_\nu G^{sc}(x_l^{(j)};x_k^s)$ and their gradients w.r.t $x^s$ by 
       \eqref{eq:disG};
      \State Obtain $G(x_l^{(j)};x_k^s) = \Phi(x_l^{(j)};x_k^s) + G^{sc}(x_l^{(j)};x_k^s)$ and  $\nabla_{x^s} G(x_l^{(j)};x_k^s) = \nabla_{x^s} \Phi(x_l^{(j)};x_k^s) + \nabla_{x^s} G^{sc}(x_l^{(j)};x_k^s)$;
       \State Evaluate $G(x;x_k^s)$ and $\nabla_{x^s} G(x;x_k^s)$ via \eqref{eq:appG2} for any $x\in \Omega_{p,H}^j$; 
      \State Evaluate $G(x;x_k^s)$ and $\nabla_{x^s} G(x;x_k^s)$ by \eqref{raylei} and then \eqref{eq:disG} for any $x\in \Omega^j_{p}\backslash\overline{\Omega_{p,H}^j}$. 
  \end{algorithmic}
\end{algorithm}

Next, we use the quasi-periodicity of $G_{\alpha}^{pq}$ to compute the
background Green's function $G$ in any other cell $\Omega_{p}^{j}$ for $j\neq
0$ and for the same set of sources $\{x_k^s\}\subset\Omega_{p,H}^0$. As in Algorithm \ref{alg:evalG}, we first evaluate $G$ on the
boundary $\partial\Omega_{p,H}^j$. To discretize $\partial\Omega_{p,H}^j$, we
horizontally translate the $(M+N)$ mesh points
$\{x_l\}_{l=1}^{M+N}\subset\partial \Omega_{p,H}^0$ by $j$ periods, denoted by $\{x_l^{(j)}\}_{l=1}^{M+N}$.
Then, ${\bm G}_\alpha^{pq,(j)}$ and ${\bm \phi}_\alpha^{pq,(j)}$, the vectors of
$G_\alpha^{pq}$ and $\phi_\alpha^{pq}$ at $x_l^{(j)}$, satisfy 
\[ 
{\bm G}_\alpha^{pq,(j)} = e^{\bi 2\pi
\alpha j }{\bm G}_\alpha^{pq},\quad  {\bm \phi}_\alpha^{pq,(j)} = e^{\bi 2\pi
\alpha j }{\bm \phi}_\alpha^{pq}. 
\]
Applying the inverse Floquet-Bloch transform
\eqref{eq:disG}, we obtain $G^{pq}$ and hence $G$ on $\partial\Omega_{p,H}^j$.
Then, similar to steps 6 and 7 in Algorithm~\ref{alg:evalG}, we can obtain $G$ and its
derivatives in $\Omega_p^j$. The above procedure is summarized in
Algorithm~\ref{alg:evalG2}.

\section{The original scattering problem}
With the background Green's function $G$ available, we are ready to solve the original scattering problem \eqref{eq:bc}-\eqref{eq:2} now. First, we construct a TBC in the perturbed region $\Omega_{np}^0$ to truncate the unbounded domain $\Omega_{np}$. 

\subsection{Transparent boundary condition}

\begin{figure}[ht!]
  \centering
  \includegraphics[width=0.9\textwidth]{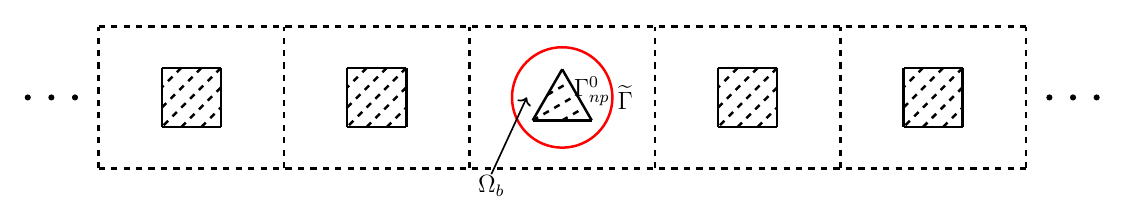}
  \caption{The artificial curve $\tilde{\Gamma}$ and the perturbed region $\Omega_b$. }
  \label{fig:5}
\end{figure}

Let $\tilde{\Gamma}$ be an artificial curve that encloses both the perturbation and the original unperturbed obstacle in the same cell. For simplicity, we choose $\tilde{\Gamma}$ to be analytic and sufficiently away from the obstacles, as indicated in Figure~\ref{fig:5}.
As we have assumed that both $u^{\rm og}$ and $G$ satisfy the SRC condition
\eqref{con:src}, we may follow the arguments in \cite[Therorem 2.5]{colkre13} to
derive the following BIE for $u^{\rm og}$ on $\tilde{\Gamma}$,
\begin{equation}
\label{eq:tbcbie}
    (\tilde{\mathcal{K}}-\mathcal{I})u^{\rm og}=\tilde{\mathcal{S}}\partial_{\nu}u^{\rm og},
\end{equation}
where the single-layer and double-layer operators $\tilde{\cal S}$ and $\tilde{\cal K}$ 
are now defined with the background Green's function $G$ in place of the fundamental solution $\Phi$,
\begin{align}
  \label{tbc:o1}
  \tilde{\mathcal{S}}[\phi](x)&=2\int_{\tilde{\Gamma}}G(x;y)\phi(y)ds(y),\\
  \label{tbc:o2}
  \tilde{\mathcal{K}}[\phi](x)&=2\text{P.V.}\int_{\tilde{\Gamma}}\partial_{\nu(y)}G(x;y)\phi(y)ds(y). 
\end{align}
The new BIE \eqref{eq:tbcbie} serve as a TBC to truncate the unbounded domain,
and can be discretized by using the same Nystr\"om method in section 3.3, but
with a rather simple procedure. We note that the adjoint double-layer and
hypersingular operators can also be incorporated here to possibly construct a
better-conditioned TBC numerically based on direct or indirect formulations
\cite{colkre13}.

Take $\tilde{\mathcal{S}}$ as an example. It satisfies
\begin{align}
  \label{tbc:o3}
  \tilde{\mathcal{S}}[\phi](x)&=2\int_{\tilde{\Gamma}}\Phi(x;y)\phi(y)ds(y)+2\int_{\tilde{\Gamma}}G^{sc}(x;y)\phi(y)ds(y).
\end{align}
As $\tilde{\Gamma}$ is smooth, we can apply the same kernel splitting
method in section 3.3 to discretize the first term without using any scaling
function. To discretize the second term, we directly apply the trapezoidal rule
as the integrand is analytic. One similarly discretizes $\tilde{\cal K}$.

Suppose $\tilde{\Gamma}$ is discretized by $N_p>0$ grid points, forming the set
$\tilde{\Gamma}_{num}$. The above discretizing procedure requires evaluating the
analytic part $G^{sc}(x;y)$ and its gradient $\nabla_y G^{sc}(x;y)$ for any
$x,y\in\tilde{\Gamma}_{num}$. This can be done by taking $N_o=N_s=N_p$ and choosing
$x_j^o,x_k^s\in \tilde{\Gamma}_{num}$ in Algorithm 3.1. Thus, we obtain two
$N_p\times N_p$ matrices $\tilde{\textbf{S}}$ and $\tilde{\textbf{K}}$
approximating the two operators $\tilde{\mathcal{S}}$ and
$\tilde{\mathcal{K}}$, respectively. Consequently, \eqref{eq:tbcbie} gives
rise to an $N_p\times N_p$ matrix
$\tilde{\textbf{N}}=(\tilde{\textbf{K}}-\textbf{I})^{-1}\tilde{\textbf{S}}$
such that 
\begin{equation}
\label{eq:tbcdis}
    {\bm u}^{\rm og} \approx \tilde{\bf N}{\bm \phi}^{\rm og},
\end{equation}
where ${\bm u}^{\rm og}$ and ${\bm \phi}^{\rm og}$ are column vectors of $u^{\rm
og}$ and $\partial_{\nu} u^{\rm og}$ at the $N_p$ grid points. Thus,
$\tilde{\bf N}$ can be regarded as an approximation of the NtD operator
$\tilde{\cal N}=(\tilde{\cal K}-{\cal I})^{-1}\tilde{\cal S}$. Based
on the NtD operator $\tilde{\cal N}$ and its approximation $\tilde{\bf
N}$, we numerically solve the original problem for the two incidences in the
following two sections. 

\subsection{Wavefield in the perturbed region}
As shown in Figure~\ref{fig:5}, let $\Gamma_{np}^0$ be
the boundary of the deformed obstacle, and $\Omega_b$ be the domain bounded by
$\tilde{\Gamma}$ and $\Gamma_{np}$. We compute $u^{\rm tot}$ in the perturbed region $\Omega_b$
in this subsection.

We study the cylindrical incidence first. Suppose the source point
$x^*\in\Omega_b$; the case $x^*\notin\Omega_b$ will be discussed in
Remark~\ref{rem:srcext}. As $u^{\rm og}(x) = u^{\text{tot}}(x;x^*)$ satisfies
the SRC condition \eqref{con:src}, equation \eqref{eq:tbcdis} implies
\begin{align}
  \label{ccc:2}
  \bm{u}^{\text{tot}}_1=\tilde{\textbf{N}}{\bm \phi}^{\text{tot}}_{1},
\end{align}
or equivalently 
\begin{align}
  \label{ccc:3} \bm{u}^{sc}_1+\bm{u}^{\text{inc}}_1=\tilde{\textbf{N}}(\bm{\phi}^{\text{s}}_{1}+\bm{\phi}^{\text{inc}}_{1}),
\end{align}
where ${\bm u}^{sc}_1$ denotes the vector of $u^{sc}$ at the $N_p$ grid points on
$\tilde{\Gamma}$, etc.  In $\Omega_b$, $u^{sc}=u^{\rm tot}-u^{\rm inc}$ satisfies the homogeneous Helmholtz
equation. Applying the same BIE approach presented in section 3.3,
we obtain an NtD matrix $\widehat{\bf N}$ that satisfies 
\begin{align}
  \label{ccc:1}
  \left.\left[\begin{array}{c}{\bm{u}^{sc}_{1}}\\{\bm{u}^{sc}_{2}}\\\end{array}\right.\right]=
  \widehat{\textbf{N}}\left[\begin{array}{c}{{\bm \phi}^{sc}_{1}}\\{{{\bm \phi}}^{sc}_{2}}\\\end{array}\right], 
\end{align}
where $\bm{u}^{sc}_{j}$ and ${\bm\phi}^{sc}_{j}$ represent column
vectors of $u^{sc}$ and $\partial_{\nu}u^{sc}$  at the grid points
on $\tilde{\Gamma}$ and $\Gamma^0_{np}$, respectively. Note that we have used the same set of
grid points as in \eqref{eq:tbcdis} and ${\bm \phi}^{sc}_2$ should
approximate a scaled normal derivative, the product of
$\partial_{\nu}u^{sc}$ and a scaling function analogous to $w$ in
\eqref{cov:1}, if $\Gamma^0_{np}$ contains corners. 
Besides, the boundary condition on $\Gamma^0_{np}$ implies
\begin{align}
  \label{ccc:4}
  \bm{u}^{\text{inc}}_2 + \bm{u}^{sc}_2=0. 
\end{align}
Solving the linear system \eqref{ccc:3}-\eqref{ccc:4}, we obtain the unknowns
$\bm{u}^{sc}_j$ and $\bm{\phi}^{sc}_j$.  Based on Green's
representation formula, analogous to \eqref{rp:1} but with $\partial\Omega_{p,H}^0$ replaced
by $\partial\Omega_b$, the boundary of $\Omega_b$, we can apply the trapezoidal rule to obtain
$u^{sc}(x)$ for any $x\in\Omega_b$ by directly using the grid points on
$\partial\Omega_b$. Then, the total wavefield $u^{\rm tot}(x;x^s) =
u^{sc}(x;x^s) + \Phi(x;x^s)$ becomes available in $\Omega_b$. 

For the plane-wave incidence $u^{\text{inc}}=e^{\bi k(\cos\theta
x_1-\sin\theta x_2)}$, we follow closely the
procedures in sections 3.2 and 3.3 to compute the quasi-periodic reference
wavefield $u^{\rm tot}_{\rm ref}$; see also \cite{lu2012high}. Then,
\eqref{ccc:3} should be changed to
\begin{equation}
  \label{ccc:5} 
  \bm{u}^{\text{tot}}_1-\bm{u}^{\text{tot}}_{\rm ref,1}=\tilde{\textbf{N}}(\bm{\phi}^{\text{tot}}_{1}-\bm{\phi}^{\text{tot}}_{\rm ref,1}),
\end{equation}
where ${\bm u}^{\rm tot}_{\rm ref,1}$ denotes the vector of $u^{\rm tot}_{\rm ref}$ at the $N_p$ grid points on $\tilde{\Gamma}$. Equation \eqref{ccc:1} still holds but with $u^{sc}$ replaced by $u^{\rm tot}$. Equation \eqref{ccc:4} now becomes ${\bm u}^{\rm tot}_2=0$. 
Solving the new linear system, we obtain ${\bm u}^{\rm tot}_j$ and ${\bm \phi}^{\rm tot}_j$ on $\partial\Omega_b$. Green's representation formula then applies to get $u^{\rm tot}$ in $\Omega_b$.

The above procedure is summarized in Algorithm~\ref{alg:utot:0}. 
\begin{algorithm}[t]
  \caption{Computing the wavefield $u^{\rm tot}$ in $\Omega_b$.} 
\label{alg:utot:0}
  \begin{algorithmic}[1]
    \Require
    $N_p$ discretization points $\{x_k^s\}_{k=1}^{N_p}$ on $\tilde{\Gamma}$; 
\State Compute $G^{sc} (x_{k'}^{s'};x_k^s)$ and $\nabla_{x^s} G^{sc}(x_{k'}^{s'};x_k^s)$ by Steps 1 through 5 of Algorithm \ref{alg:evalG} for any $1\leq k',k\leq N_p$;  
\State Approximate the NtD operator $\tilde{\cal N}$ by the $N_p\times N_p$ matrix $\tilde{\bf N}$ for the exterior domain $\Omega_{np}\backslash\overline{\Omega_b}=\Omega_{p}\backslash\overline{\Omega_b}$, and establish the TBC equation \eqref{eq:tbcdis};
\State Use the BIE approach in section 3.3 to obtain the NtD matrix $\hat{\bf N}$ in \eqref{ccc:1} for the interior domain $\Omega_b$; 
\State Solve equations \eqref{ccc:3}-\eqref{ccc:4} for ${\bm u}_j^{\rm sc}$ and ${\bm \phi}_j^{\rm sc}$ for $j=1,2$ on $\tilde{\Gamma}\cup \Gamma^0_{np}$; 
\State For any $x\in \Omega_b$, evaluate $u^{\rm sc}(x)$ by using the trapezoidal rule to discretize
\begin{align*}
u^{\rm tot}(x)=&u^{\rm inc}(x) + \int_{\partial \Omega_b}[\partial_{\nu(y)} \Phi(x;y)u^{\rm sc}(y)-\Phi(x;y)\partial_{\nu} u^{\rm sc}(y)]ds(y),
\end{align*}
in terms of ${\bm u}_j^{\rm sc}$ and ${\bm \phi}_j^{\rm sc}$ for $j=1,2$;
  \end{algorithmic}
\end{algorithm}



\subsection{Wavefield in the unperturbed region}
In this section, we shall make use of Algorithms~\ref{alg:evalG} and
\ref{alg:evalG2} to compute the wavefield outside $\Omega_b$. To simplify
the presentation, we discuss the cylindrical incidence only; the case for the
plane-wave incidence can be analyzed similarly.

\begin{figure}[H]
  \centering
  \includegraphics[width=0.9\textwidth]{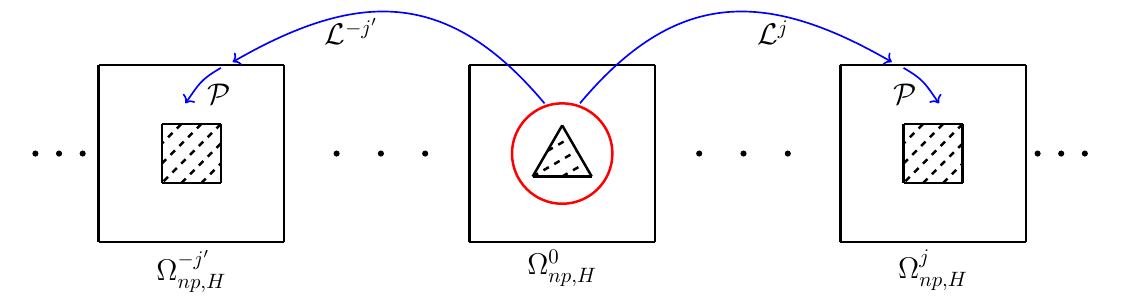}
  \caption{The leap procedures ${\cal L}^{j}$ and  ${\cal L}^{-j'}$, and the pullback procedure
  ${\cal P}$.}
  \label{fig:7}
\end{figure}

The algorithm can basically be divided into two
procedures: (1) Leap; (2) Pullback, as illustrated in Figure~\ref{fig:7}. In the leap procedure, we make use of the
total wavefield $u^{\rm tot}$ on the artificial boundary $\tilde{\Gamma}$ 
to get $u^{\rm tot}$ on the boundary $\partial\Omega^{j}_{np,H}$ of the $j$-th cell,
called the leap procedure ${\cal L}_j$.  In the pullback procedure, we make use
of $u^{\rm tot}$ on $\partial\Omega^{j}_{np,H}$ to compute $u^{\rm tot}$ in the
$j$-th unit cell $\Omega^{j}_{np}$, called the pullback procedure ${\cal P}$. 
We discuss $j\neq 0$ first so that $\Omega_{np}^j = \Omega_{p}^j$. To evaluate
$u^{\rm tot}$ in $\Omega^j_{p}$, we need to compute $u^{\rm tot}(x)$ and the
normal derivative $\partial_\nu u^{\rm tot}(x)$ for $x\in \partial
\Omega^j_{p,H}$ first; note that $u^{\rm tot}=0$ on the obstacle boundary
$\Gamma_{p}^j$. As we have obtained $u^{\rm tot}$ and $\partial_\nu u^{\rm tot}$
on $\tilde{\Gamma}$, we directly adopt the following Green's representation
formula 
\begin{align} 
\label{eq:u:gamma1} 
u^{\rm tot}(x)=\int_{\tilde{\Gamma}}[\partial_{\nu(y)} G(x;y)u^{\rm tot}(y)-G(x;y)\partial_{\nu(y)} u^{\rm tot}(y)]ds(y), 
\end{align} 
to evaluate $u^{\rm tot}(x)$ for $x\in \partial \Omega^j_{p,H}$. In doing so, we use the
trapezoidal rule to discretize the above integral by using
$\tilde{\Gamma}_{num}$ as the set of grid points on $\tilde{\Gamma}$ and then
apply Algorithm 3.1 to evaluate $\partial_{\nu(y)} G(x;y)$ and $G(x;y)$ for
$y\in \tilde{\Gamma}_{num}$ and $x\in \partial \Omega^j_{p,H}$. We can directly take the
gradient of \eqref{eq:u:gamma1} w.r.t. $x$ to evaluate $\nabla u^{\rm tot}$ and
hence $\partial_\nu u^{\rm tot}$ on $\partial\Omega^1_{p,H}$ as well.  Next, Green's
representation formula, using $\Phi$ instead of $G$, applies to get $u^{\rm
tot}$ everywhere in $\Omega^j_{p,H}$. As for $u^{\rm tot}(x)$ for any $x\in
\Omega^j_{p}\cap\{x:|x_2|\geq H\}$, we can use step 6 of
Algorithm~\ref{alg:evalG2} to evaluate $G(x;y)$ and its derivatives for
$y\in\tilde{\Gamma}_{num}$, and then use \eqref{eq:u:gamma1} to compute
$u^{\rm tot}(x)$. As for $j=0$, we only need to replace $\partial
\Omega_{p,H}^j$ above by the rectangular boundary
$\partial((-\pi,\pi)\times(-H,H))$.  We put all details together in Algorithm
\ref{alg:utot:j} below.
\begin{algorithm}[H]
  \caption{Computing the wavefield in $\Omega^j_{np}\backslash\overline{\Omega_b}$ for any $j\in\mathbb{Z}$.} 
\label{alg:utot:j}
  \begin{algorithmic}[1]
    \Require
    $N$ mesh points $\{x_l^{(j)}\}_{l=1}^{N}$ on $\partial\Omega^{j}_{np,H}$; $N_p$ mesh points $\{x_k^s\}_{k=1}^{N_p}$ on $\tilde{\Gamma}$; $u^{\rm tot}(x_k^s)$ and $\partial_\nu u^{\rm tot}(x_k^s)$ on $\tilde{\Gamma}$ computed in section 4.2.
\State Compute $G(x_l^{(j)};x_k^s)$ and $\nabla_{x^s} G(x_l^{(j)};x_k^s)$ by Steps 1 through 4 of Algorithm \ref{alg:evalG2};  
\State Leap ${\cal L}_j$: evaluate $u^{\rm tot}(x_l^{(j)})$ and $\partial_\nu u^{\rm tot}(x_l^{(j)})$ on $\partial\Omega_{np,H}^j\backslash\overline{\Omega_b}$ by using the trapezoidal rule, with the set of nodes $\tilde{\Gamma}_{\rm num}$, to discretize
\begin{align*}
\label{rp:4.4_1}
\partial_{x_1}^{j_1}\partial_{x_2}^{j_2}u^{\rm tot}(x)=&\partial_{x_1}^{j_1}\partial_{x_2}^{j_2}\int_{\tilde{\Gamma}}[\partial_{\nu(y)} G(x;y)u^{\rm tot}(y)-G(x;y)\partial_{\nu} u^{\rm tot}(y)]ds(y), 
\end{align*}
for $0\leq j_1+j_2\leq 1$.
\State Pullback ${\cal P}$: for any $x\in \Omega_{np,H}^j\backslash\overline{\Omega_b}$, 
evaluate $u^{\rm tot}(x)$ by using the trapezoidal rule to discretize
\begin{align*}
u^{\rm tot}(x)=&\int_{\partial(\Omega^{j}_{np,H}\backslash\overline{\Omega_b})}[\partial_{\nu(y)}\Phi(x;y)u^{\rm tot}(y)-\Phi(x;y)\partial_{\nu} u^{\rm tot}(y)]ds(y).
\end{align*}
\State For any $x\in \Omega^j_{np}\backslash\overline{\Omega^j_{np,H}}$, compute $G(x;x_k^s)$ and $\nabla_{x^s} G(x;x_k^s)$ by Step 6 of Algorithm \ref{alg:evalG2}, and evaluate $u^{\rm tot}(x)$ by using the trapezoidal rule, with the set of nodes $\tilde{\Gamma}_{\rm num}$, to discretize
\begin{align*}
u^{\rm tot}(x)=&\int_{\tilde{\Gamma}}[\partial_{\nu(y)} G(x;y)u^{\rm tot}(y)-G(x;y)\partial_{\nu} u^{\rm tot}(y)]ds(y).
\end{align*}
  \end{algorithmic}
\end{algorithm}

Throughout Algorithms~\ref{alg:utot:0} and \ref{alg:utot:j}, we observe that the
background Green's function $G(x;y)$ is evaluated for $x\in \tilde{\Gamma}\cup
(\cup_{j\in\mathbb{Z}} \partial\Omega_{p,H}^j)$ and $y\in \tilde{\Gamma}$. Owing
to the IFB transform \eqref{eq:ifbG} and  the quasi-periodicity of
$G^{qp}_\alpha$, it is enough to evaluate $G^{qp}_\alpha(x;y)$ for $x\in
\tilde{\Gamma}\cup \partial\Omega_{p,H}^0$, $y\in \tilde{\Gamma}$, and
$\alpha\in \{\alpha_j\}_{j=1}^{n}$ only. This is the fundamental feature that
accounts for the efficiency of the proposed solver! 

We make two remarks to conclude this section.
\begin{myremark}
\label{rem:genptb}
We have assumed so far that only one cell is perturbed. In
fact, the current approach can handle more general situations.  For example,
when there are two perturbed cells that are far apart, the scattering problem
becomes numerically challenging to solve if we use a single large curve
$\tilde{\Gamma}$ to enclose both perturbed regions. Instead, we can use two
separate smaller curves, each enclosing the perturbed region within its
respective cell, as illustrated in Figure~\ref{fig:8}. 
\begin{figure}[ht!]
  \centering
  \includegraphics[width=0.9\textwidth]{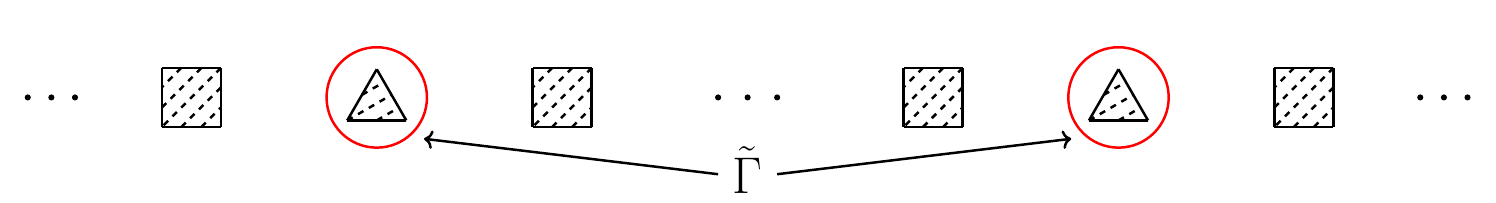}
  \caption{Two far apart, perturbed unit cells.}
  \label{fig:8}
\end{figure}
By setting $\tilde{\Gamma}$ as the union of the two separate curves, we can
follow exactly the same procedures to compute the wavefield in all cells. This
should be compared with \cite{yu2022pml, joly2006exact}. The numerical solvers
therein are based on two marching operators to truncate the periodic direction, and
have to discretize all cells between the two perturbed cells before
computing the wavefield everywhere. 
\end{myremark}
\begin{myremark}
\label{rem:srcext}
If $x^*\notin\Omega_b$, certainly, we can follow the same approach as in
Remark~\ref{rem:genptb}, using an extra closed curve to enclose the source point
$x^*$. However,  discretizing this extra curve increases the computational
burden. To resolve this issue,  we redefine $u^{\rm og}(x;x^*)= u^{\rm
tot}(x;x^*) - G(x;x^*)$ for $x\notin\Omega_b$, where $G(x;x^*)$ is evaluated by
Algorithms~\ref{alg:evalG} and \ref{alg:evalG2},  and the current approach can
be easily adapted to compute $u^{\rm tot}$. 
\end{myremark}

\section{Numerical examples} 
In this section, we carry out several numerical experiments to validate the
proposed numerical solver. To illustrate the accuracy of the method, we first
compute the wavefield $u^{\rm tot}$ on the artificial curve $\tilde{\Gamma}$ for
some sufficiently large $n$, $M$ and $N$, where we recall that $n$ is the number
of nodes in the discretized IFB transform \eqref{eq:disG}, $M$ is the total
number of grid points on the obstacle boundary $\Gamma_{np}^j$ in the $j$-th
unit cell, and $N$ is the total number of grid points on the rectangular
boundary $\partial\left[(2\pi j -\pi,2\pi j+\pi)\times(-H,H)\right]$.
Then, we compare it with the less accurate numerical solutions for some smaller
values of $n$, $M$ or $N$,  at a common set of grid points on $\tilde{\Gamma}$
by checking the relative error
\begin{align*}
  \text{E}_{\rm rel}=\frac{||{\bf u}^{\rm tot}_{\rm num}-{\bf u}^{\rm tot}_{\rm ref}||_{\infty}}{||{\bf u}^{\rm tot}_{\rm ref}||_{\infty}},
\end{align*}
where ${\bf u}^{\rm tot}_{\rm ref}$ denotes the vector of the reference solution
on $\tilde{\Gamma}$ and  ${\bf u}^{\rm tot}_{\rm num}$ denotes the vector of a
numerical solution to be compared.

We shall demonstrate the performance of the current solver for two different
wavenumbers $k=1.25$ and $k=1.5$, and compare it with our recently proposed
PML-BIE solver \cite{yu2022pml}, as \cite{zhang2022fast} hints a poorer
performance of the PML truncation for half-integers. In all examples, we assume
the period is $2\pi$, take $p=6$ in \eqref{cov:1}, and use $\tilde{N}_p = 300$
points to discretize the artificial curve $\tilde{\Gamma}$, which ensures that
the discretization error for the TBC is sufficiently small.

\begin{figure}[htbp]
  \centering
  \includegraphics[width=0.9\textwidth]{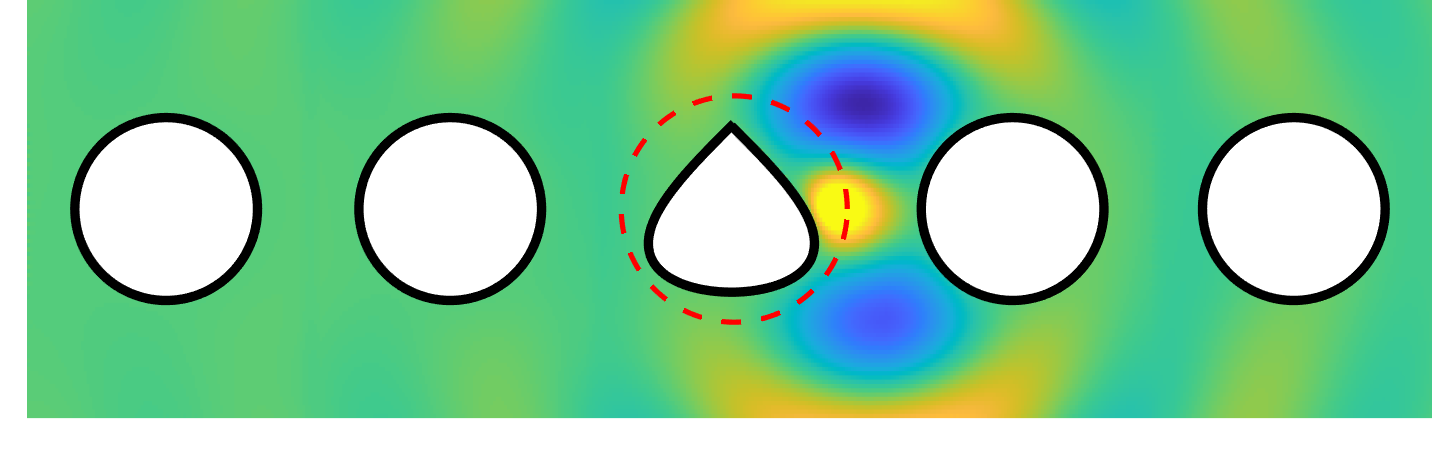}
  \caption{Example 1: Real part of $u^{\rm tot}$ for $k=1.25$.}
  \label{fig:8.1}
 \end{figure}

\noindent{\bf Example 1.} In this example, the obstacles are disks of radius $2$
and one of them is replaced by a drop-shaped obstacle, with the boundary
$\{x=(1.8\sin(t),3.6\sin(\frac{1}{2}t)+1.8),0\leq t\leq 2\pi\}$. The boundaries
of these obstacles are illustrated by the solid lines in Figure~\ref{fig:8.1}.
We compute the total wavefield $u^{\rm tot}$ excited by a source at $x^*=(2.2, 0)$. 

\begin{figure}[htb]
    \centering
    (a)\includegraphics[width=0.4\textwidth]{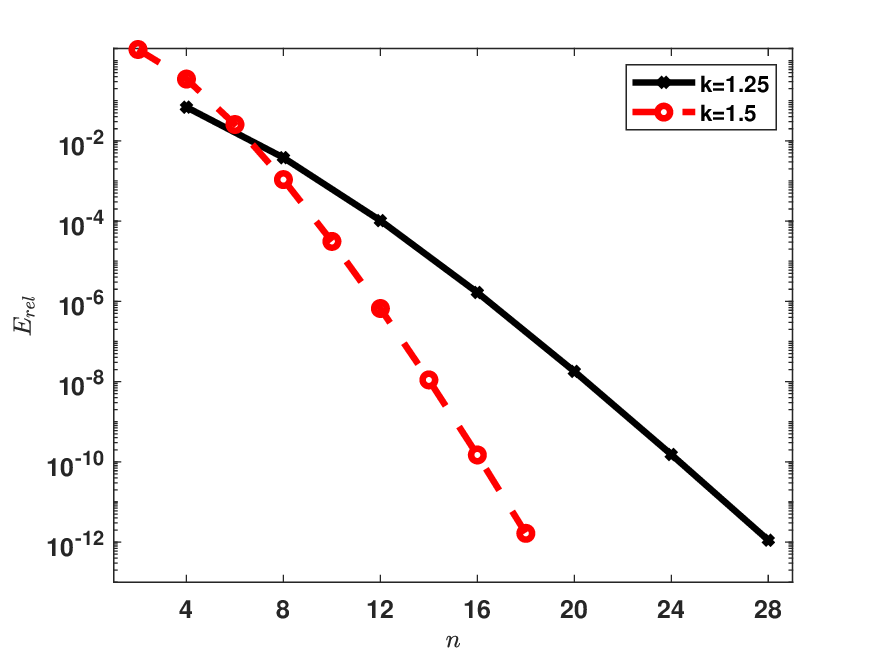}
    (b)\includegraphics[width=0.4\textwidth]{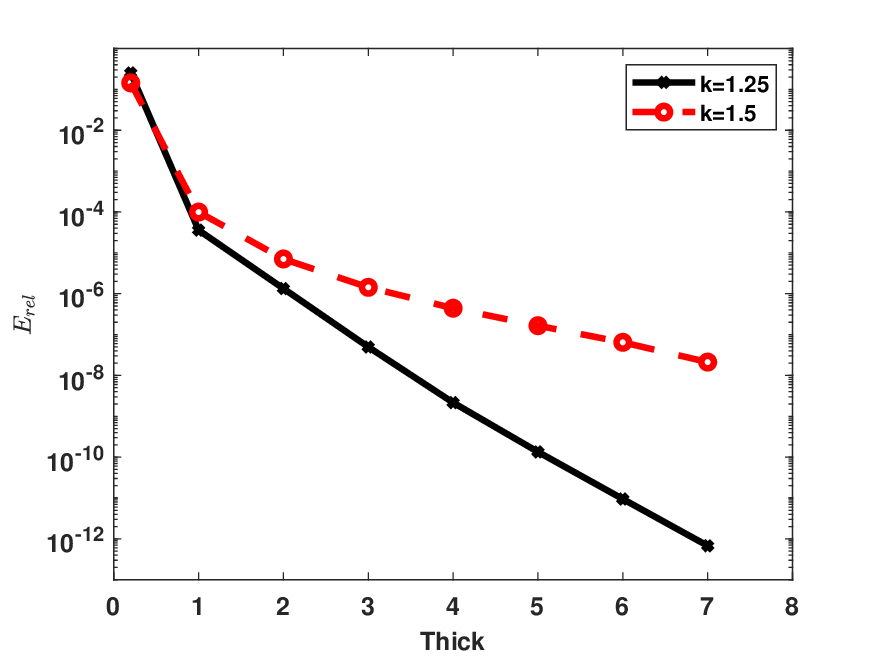}
    \caption{Example 1: (a) $E_{\rm rel}$ against $n$ for $k=1.25,1.5$ based on the current method; (b) $E_{\rm rel}$ against $n$ for $k=1.25,1.5$ based on the PML-BIE method \cite{yu2022pml}.}
    \label{fig:7.1.2}
\end{figure}

As illustrated by the dashed lines in Figure~\ref{fig:8.1}, the artificial curve
$\tilde{\Gamma}$ is chosen as a circle of radius $2.5$ centered at $(0,0)$ that encloses the
drop-shaped obstacle.  By choosing $H=\pi$, $M=150$, $N=600$, we compute the
total wavefield $u^{tot}$ for the wavenumber $k=1.25$ with $n=32$ in
\eqref{eq:disG} and for the wavenumber $k=1.5$ with $n=18$, and obtain two
reference solutions; note that for the half-integer wavenumber, the IFB
transform involves a single integral so that a fewer number of discretization
points are needed. The field pattern of the reference solution for $k=1.25$ in
five consecutive cells is shown in Figure~\ref{fig:8.1}. We evaluate $E_{\rm
rel}$ by comparing the numerical solutions with the reference solutions at $300$
uniformly spaced grid points on $\tilde{\Gamma}$ for different values of $n$.
The convergence curves for the two different wavenumbers are shown in
Figure~\ref{fig:7.1.2}(a).

To apply the PML-BIE solver~\cite{yu2022pml}, we place two PMLs above $x_2=H$
and below $x_2=-H$, and evaluate the total wavefield $u^{\rm tot}$ on
$\tilde{\Gamma}$ using exactly the same number of points on the boundary of each
unit cell in the physical domain. Based on the computed two reference
solutions, we evaluate $E_{\rm rel}$ for the PML solutions for different values
of PML thickness. The corresponding convergence curves for the two different
wavenumbers are shown in Figure~\ref{fig:7.1.2}(b). 

\begin{figure}[htb]
    \centering
    (a)\includegraphics[width=0.45\textwidth]{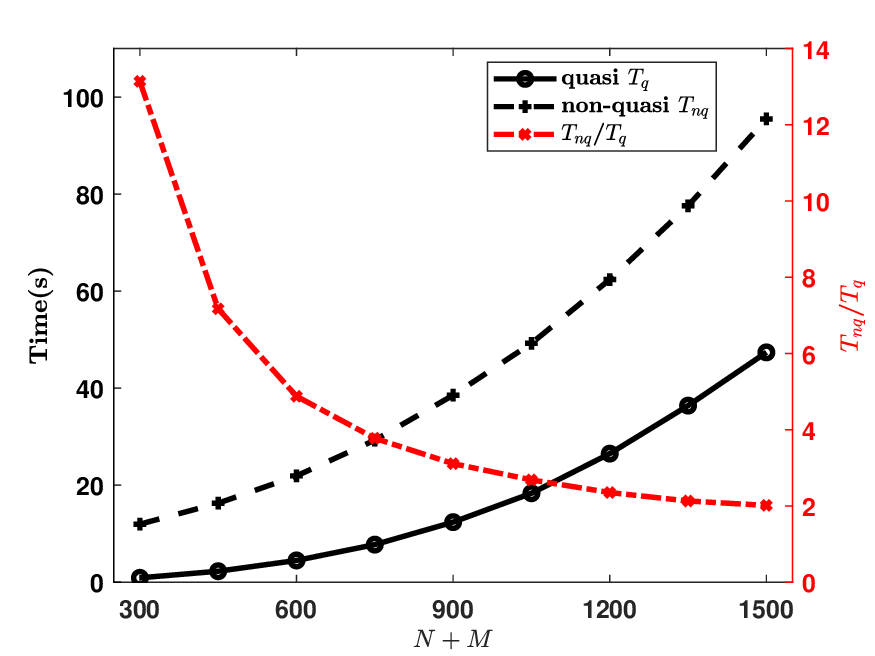}
    (b)\includegraphics[width=0.45\textwidth]{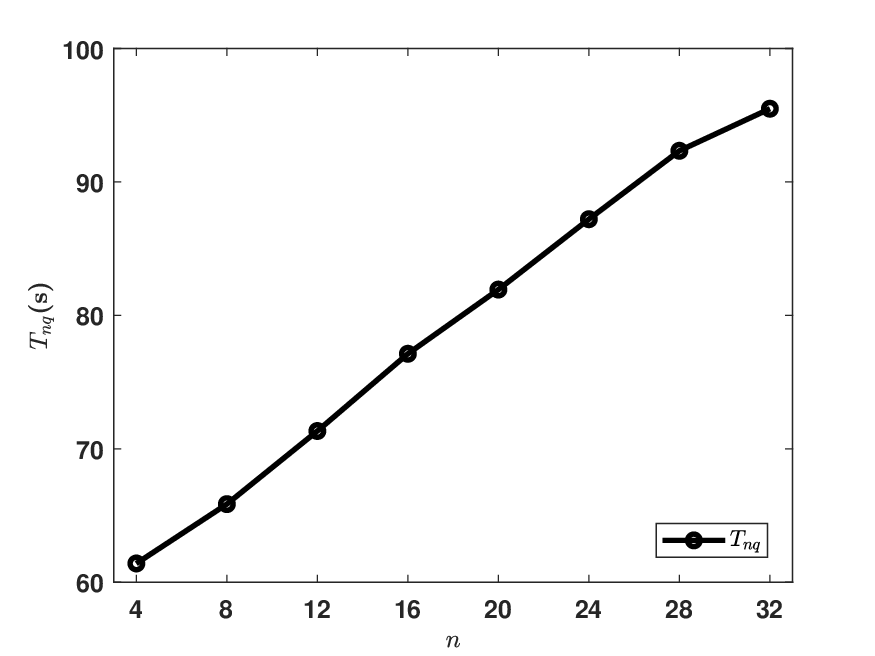}
    \caption{Example 1: (a) Running times $T_{nq}$ and $T_q$ for the non-quasi-periodic problem and a quasi periodic problem against $N+M$; (b) $T_{nq}$ against $n$. }
    \label{fig:7.1.1}
\end{figure}

Compare Figures~\ref{fig:7.1.2} (a) and (b). When $k$ is not a half-integer,
both methods can provide roughly 12 accurate digits. However, we find that
the current method takes only 40 seconds, while the PML-BIE method takes more than 4
minutes.  When $k$ is a half-integer, the PML-BIE method provides only 7
accurate digits for the same thickness and its convergence rate
downgrades significantly. Nevertheless, the current method can still provide
$12$ accurate digits.
In Figure~\ref{fig:7.1.1}(a), we also show the running times $T_{nq}$ for solving
the non-quasi-periodic problem with $n=32$ and $N_p=300$ and $T_q$ for solving one single quasi-periodic problem, for different values of
$M+N$. It can be seen that the running time of the non-quasi periodic problem is
only about twice that of one quasi periodic problem  for sufficiently large $M+N$,
demonstrating the efficiency of the current solver. 
Figure~\ref{fig:7.1.1}(b) shows the variation of $T_{nq}$ for different values of $n$ with $N+M=1500$ and $N_p=300$. It can be seen that $T_{nq}$ does not increase in direct proportion to $n$. Instead, an eightfold increase in $n$ results in a less than twofold increase in $T_{nq}$ only. 

\begin{figure}[htbp]
  \centering
  \includegraphics[width=0.9\textwidth]{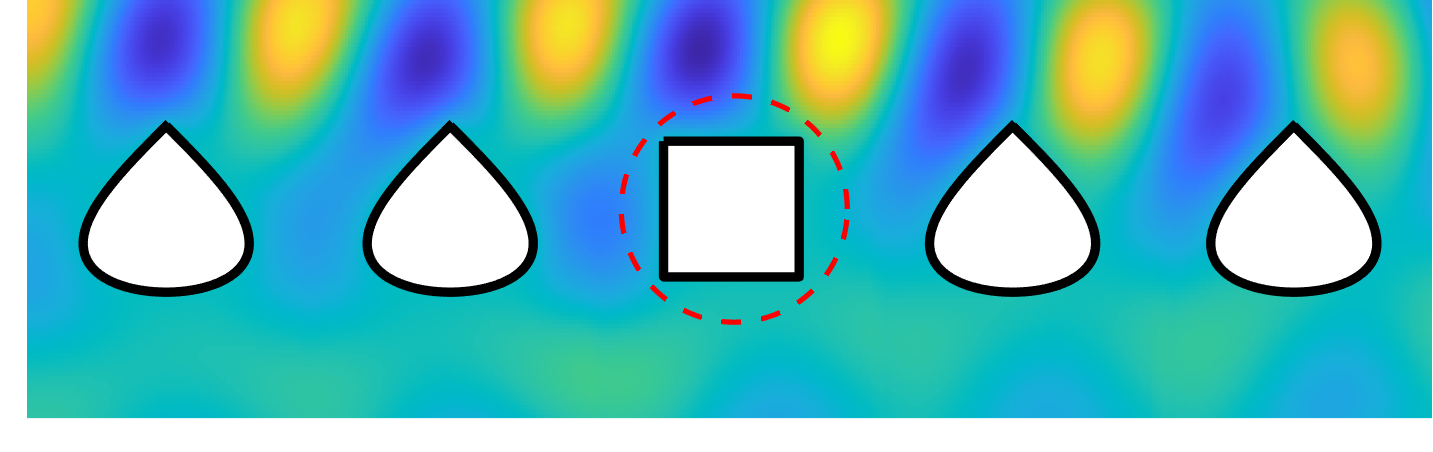}
  \caption{Example 2: Real part of $u^{\rm tot}$ for $k=1.25$.}
  \label{fig:8.5}
 \end{figure}

\noindent{\bf Example 2.} In this example, we use the drop-shaped obstacle in Example 1 to construct a $2\pi$-periodic structure and replace one of them by a square obstacle with a side length of 3, as illustrated by the solid lines in Figure~\ref{fig:8.5}.
We compute the total wavefield $u^{\rm tot}$ excited by a plane wave 
$e^{\bi k(\cos\theta x_1-\sin\theta x_2)}$ with the incident angle $\theta=\frac{\pi}{6}$. 
As illustrated by the dashed lines in Figure~\ref{fig:8.5}, the artificial curve $\tilde{\Gamma}$ is also chosen as a circle of radius $2.5$ that encloses the
square obstacle. By choosing the same parameters $H$, $M$, $N$ and $n$ as example 1, we compute reference solutions of the total wavefield $u^{tot}$ for two different wavenumbers $k=1.25$ and $k=1.5$. The field pattern of $\Re(u^{\rm
tot})$ for $k=1.25$ in five consecutive cells is shown in Figure~\ref{fig:8.5}.

\begin{figure}[htbp]
    \centering
    (a)\includegraphics[width=0.4\textwidth]{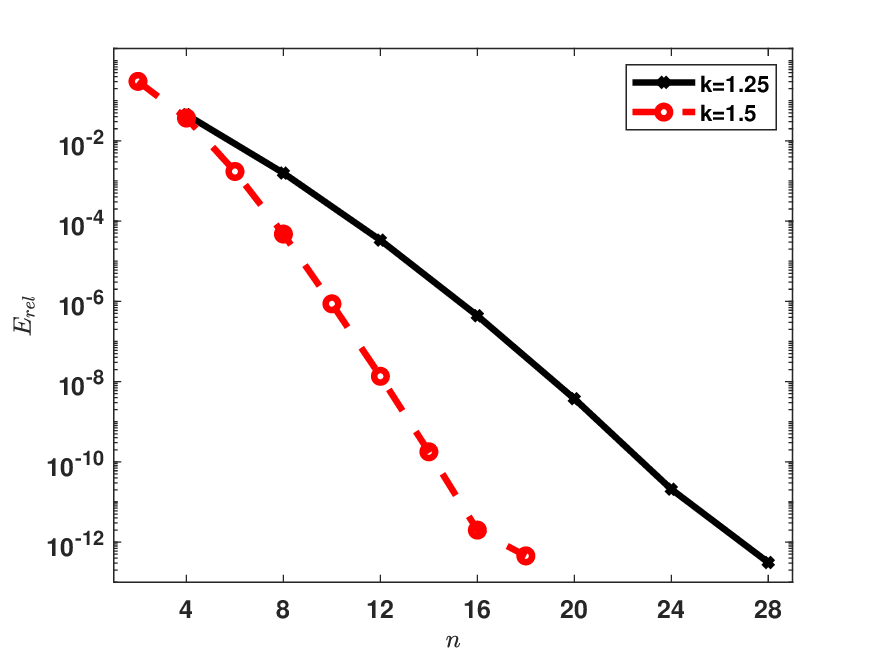}
    (b)\includegraphics[width=0.4\textwidth]{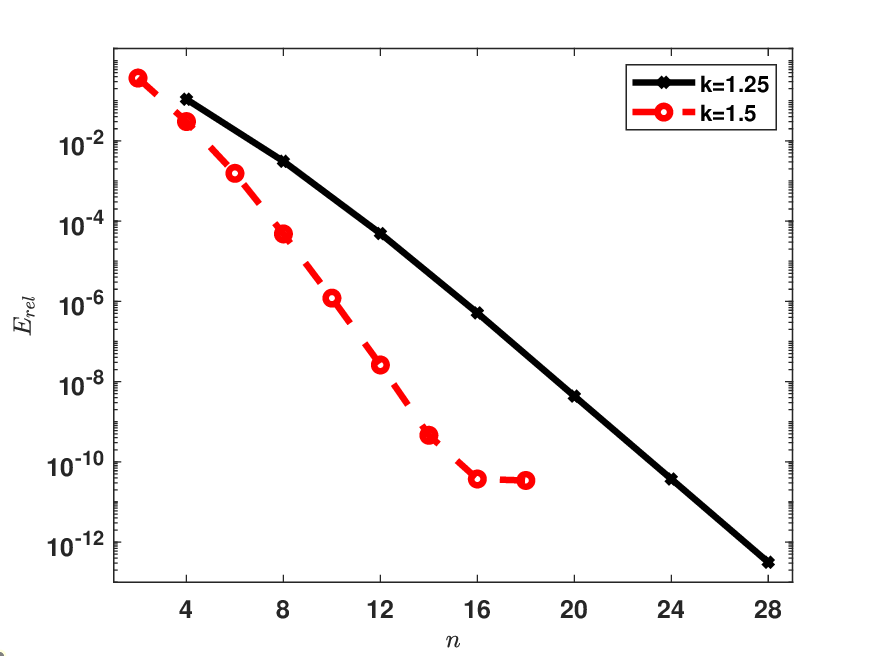}
    \caption{
    (a) $E_{\rm rel}$ against $n$ for $k=1.25,1.5$ for Example 2; (b) $E_{\rm rel}$ against $n$ for $k=1.25,1.5$ for Example 3. }
    \label{fig:7.1.3}
\end{figure}

We evaluate $E_{\rm rel}$ by comparing the numerical solutions with the reference solutions at $300$ grid points on $\tilde{\Gamma}$ for different values of $n$. The
convergence curves for the two different wavenumbers are shown in Figure~\ref{fig:7.1.3}(a).

\begin{figure}[htbp]
  \centering
  \includegraphics[width=0.9\textwidth]{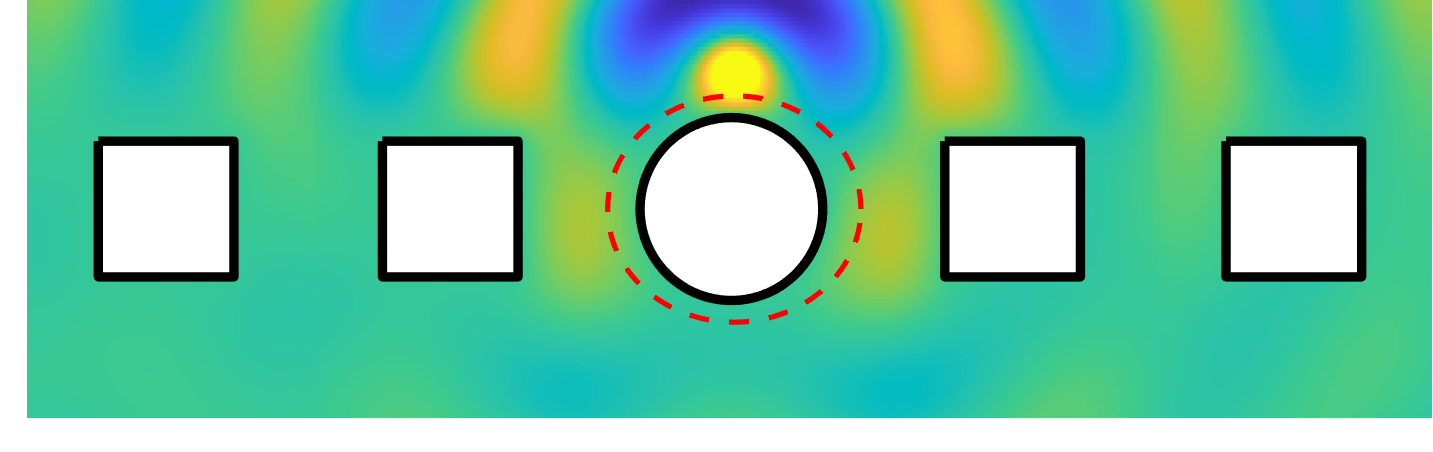}
  \caption{Example 3: Real part of $u^{\rm tot}$ for $k=1.25$.}
  \label{fig:8.3}
\end{figure}

\noindent{\bf Example 3.}
In this example, we use the square obstacle in Example 2 to produce a $2\pi$-period structure and replace one of them by a disk of radius $2$, 
as illustrated by the solid lines in Figure~\ref{fig:8.3}.
We compute the total wavefield $u^{\rm tot}$ excited by a source 
$x^*=(0,3)$. 
As illustrated by the dashed lines in Figure~\ref{fig:8.3}, the artificial curve $\tilde{\Gamma}$ is chosen as an ellipse with a major axis of $5.6$ in $x_1$-direction and a minor axis of $5$ in $x_2$-direction that encloses the disk obstacle. 
By choosing $H=\pi$, $M=600$ and $N=600$, and $n=32$ for $k=1.25$ and $n=18$ for $k=1.5$, we obtain two reference solutions. The field pattern of $\Re(u^{\rm
tot})$ for $k=1.25$ in five consecutive cells is shown in Figure~\ref{fig:8.3}.
We evaluate $E_{\rm rel}$ by comparing numerical solutions with the reference solution at $300$ grid points on $\tilde{\Gamma}$ for different values of $n$. The
convergence curves for the two different wavenumbers are shown in Figure~\ref{fig:7.1.3}(b).


\section{Conclusion}

In this paper, we have developed a new BIE method based on the FB transform for
solving an acoustic wave scattering problem in a locally perturbed periodic
structure in two dimensional spaces. Adopting the FB and IFB transforms, we developed
efficient algorithms for computing the Green's functions and its derivatives for
the background periodic structure. These Green's functions helped to construct a
BIE on an artificial curve enclosing the perturbation, serving as a TBC to
truncate the unbounded domain. The BIE was discretized by a spectral accuracy
quadrature rule. Effective algorithms based on leap and pullback procedures were
further developed to compute the total wavefield everywhere in the structure.
Numerical experiments have demonstrated the efficiency and accuracy of this new
method.

The method can be extended towards a number of directions. For periodic
obstacles that are sound-hard or penetrable, and for periodic surfaces used in
diffraction gratings, as long as there are no guided modes, this approach can be
extended to construct the TBC. In particular, when the periodic structure
contains unbounded periodic interfaces, the TBC will be defined across different
homogeneous media, making the background Green's function more difficult to evaluate
in the vicinity of the interfaces. On the other hand, if guided modes indeed exist,
the contour deformation approach for the IFB transform in
\cite{zhang2021numerical}, constructing a new path to detour the guided modes,
can be used to construct the TBC. Nevertheless, we are more interested in
dealing with the real-path integrals with integrands of poles directly
\cite{kirsch2018limiting}, as such poles cannot be detoured for locally
perturbed bi-periodic structures in three dimensions \cite{arens2024high}. We wish to tackle
these issues in the future.

\section*{Acknowledgment}
This work was done when the second author K. Shen was visiting the third author R. Zhang at TU Berlin. Shen is grateful to Prof. Zhang for her kind help and insightful supervision during the visit. W. Lu is partially supported by the National Key Research and Development Program of China (2024YFA1012600, 2023YFA1009100), China NSF grant 12174310, and a Key Project of Joint Funds for Regional Innovation and Development (U21A20425), Zhejiang Provincial NSF Distinguished Youth Program (Extended Program LRG25A010001).





\bibliographystyle{cas-model2-names}

\bibliography{cas-refs}






\end{document}